\pdfoutput=1

\newcommand{\kernel}{\widetilde{\operatorname{H}}}  % Generic Kernel-real var.
\newcommand{\kernelu}{\widetilde{\operatorname{H}}} % Generic Kernel-first var.
\newcommand{\kerneluu}{{\operatorname{H}}}          % Generic Kernel-final var.

\newcommand{\surfelu}{\widetilde{J}} % Surface Element-first var.
\newcommand{\surfeluu}{{J}}          % Surface Element-final var.

\newcommand{\densiw}{\psi}  % Modified density
\newcommand{\densicheb}{a^q}  % Chebyshev coefficients

\newcommand{\chebyw}[4]{\beta_{#2 , #3}^{\, #1 , #4}} % Order of arguments(q,n,m,\ell)

\newcommand{\field}{U}                 % Acoustic field
\newcommand{\densi}{\,\widetilde{\varphi}\,}  % Density-real var.
\newcommand{\densiu}{\widetilde{\varphi}} % Density-first var.
\newcommand{\densiuu}{{\varphi}}          % Density-final var.

\newcommand{\npts}[2]{{N_#1^{\,#2}}}

\newcommand{\wck}{w} % Change of varibles

\newcommand{\fejernodes}{x}
\newcommand{\fejerw}{w} % Fej\'er weights

\newcommand{\setclose}{\Omega^c} % Set for points close to the patch
\newcommand{\setfar}{\Omega^f}   % Set for points far to the patch
\newcommand{\setpts}{\chi}       % Set for discretization points

\newcommand{\paramorig}[1]{\widetilde{\nex}^{\,#1}}
\newcommand{\param}[1]{\nex^{\,#1}}

\newcommand{\uorig}{s} % First var-u
\newcommand{\vorig}{t} % First var-v

\newcommand{\wu}{{\eta_s^{\, q}}} % Change of variable corner-u
\newcommand{\wv}{{\eta_t^{\, q}}} % Change of variable corner-v

\newcommand{\uproj}[2]{{\overline{\uu}^{\,#1}_{#2}}}
\newcommand{\vproj}[2]{{\overline{\vv}^{\,#1}_{#2}}}

\newcommand{\wsingu}[1]{{\xi_{#1}}} % Change of variable singular
\newcommand{\wsinguweight}{{\mu}}

\newcommand{\nsing}{N_{\beta}}

\newcommand{\uu}{u} % Final var-u
\newcommand{\vv}{v} % Final var=v

\newcommand{\bdry}{\Gamma}
\newcommand{\dsurf}{\text{d} \sigma(\ney)}

\newcommand{\reals}{\mathbb{R}}

\newcommand{\spot}{\mathscr{S}}
\newcommand{\dpot}{\mathscr{D}}

\newcommand{\soper}{S}
\newcommand{\doper}{D}

\newcommand{\de}{\text{d}}
\newcommand{\p}{\partial}
\newcommand{\ney}{\boldsymbol{r}'}
\newcommand{\nex}{\boldsymbol{r}}
\newcommand{\nn}{\boldsymbol{n}}
\newcommand{\inc}{\text{inc}}
\newcommand{\scat}{\text{scat}}
\newcommand{\bigo}{\mathcal{O}}
\documentclass[11pt,letter]{article}

\usepackage{float}
\usepackage[numbers]{natbib}

\usepackage[colorlinks=true,allcolors=blue]{hyperref}
\usepackage{latexsym}
\usepackage{amsmath}
\usepackage{amsfonts}
\usepackage{amssymb}
\usepackage{amsthm}
\usepackage{graphicx}
\usepackage{empheq}
\usepackage{caption}
\usepackage{enumerate}
\usepackage{color}
\usepackage{subfig}
\usepackage{mathrsfs}

\begin{document}

\title{A Chebyshev-based rectangular-polar integral solver for scattering by
  general geometries described by non-overlapping patches}

\author{Oscar P.
  Bruno$^1\footnote{E-mail: \texttt{\{obruno,emmanuel.garza\}@caltech.edu}}$,
  Emmanuel Garza$^1$\\\\
(Preliminary Draft)\\ \\
  \small{$^1$Computing \& Mathematical Sciences, California Institute of
    Technology}}

\date{\today}

\maketitle

\begin{abstract}
  This paper introduces a high-order-accurate strategy for integration
  of singular kernels and edge-singular integral densities that appear
  in the context of boundary integral equation formulations of the
  problem of acoustic scattering. In particular, the proposed method
  is designed for use in conjunction with geometry descriptions given
  by a set of arbitrary non-overlapping logically-quadrilateral
  patches---which makes the algorithm particularly well suited for
  treatment of CAD-generated geometries. Fej\'er's first quadrature rule
  is incorporated in the algorithm, to provide a spectrally accurate
  method for evaluation of contributions from far integration regions,
  while highly-accurate precomputations of singular and near-singular
  integrals over certain ``surface patches'' together with
  two-dimensional Chebyshev transforms and suitable surface-varying
  ``rectangular-polar'' changes of variables, are used to obtain the
  contributions for singular and near-singular interactions. The
  overall integration method is then used in conjunction with the
  linear-algebra solver GMRES to produce solutions for sound-soft
  open- and closed-surface scattering obstacles, including an
  application to an aircraft described by means of a CAD
  representation. The approach is robust, fast, and highly accurate:
  use of a few points per wavelength suffices for the algorithm to
  produce far-field accuracies of a fraction of a percent, and slight
  increases in the discretization densities give rise to significant
  accuracy improvements.
\end{abstract}
%-------------------------------------------------------------------------------

%-------------------------------------------------------------------------------

\section{Introduction} \label{sec:intro} %--------------------------------------

The solution of scattering problems by means of boundary integral
representations has proven to be a game-changer when the ratio of volume to
surface scattering is large, where volumetric solvers become intractable due to
memory requirements and computational cost. At the heart of every boundary
integral equation (BIE) solver lies an integration strategy that must be able to
handle the weakly singular integrals associated with the Green-function
integral formulations of acoustic and electromagnetic scattering. Several
approaches have been proposed to deal with this difficulty, most notably those
put forward in~\citep{Kunyansky-1,Ganesh-1,Bremer-1}.

For three-dimensional scattering, which reduces to two-dimensional
integral equations along the scatterer's surface, there is no simple
quadrature rule that accurately evaluates the weakly singular
scattering kernels, which makes the three-dimensional problem
considerably more difficult than its two dimensional counterpart---for
which a closed form expression exists for the integrals of the product
of certain elementary singular kernels and complex
exponentials~\citep{Colton-1}. Therefore, a number of approaches have
been proposed to treat these integrals---which can be classified into
the Nystr\"om, collocation and Galerkin
categories~\citep{Colton-1}. Nystr\"om methods use a quadrature rule
to evaluate integrals, which leads to a system of linear equations;
the collocation approach finds a solution on a finite-dimensional
space which satisfies the continuous BIE at some collocation points
(hence the name); the Galerkin approach solves the BIE in the weak
formulation using finite-element spaces for both the solution and test
functions.

In this contribution, we present what can be considered as a hybrid
Nystr\"om-collocation method, in which the far interactions are
computed via Fej\'er's fist quadrature rule (Nystr\"om), which yields
spectrally accurate results for smooth integrands, while the integrals
involving singular and near-singular kernels are obtained by relying
on highly-accurate precomputations of the kernels times Chebyshev
polynomials (which are produced by means of a surface-varying
rectangular-polar change of variables), together with Chebyshev
expansions of the densities (collocation). The derivatives of the
rectangular-polar change of variables vanish at the kernel singularity
(the surface-varying observation point) and geometric-singularity
points, producing a ``floating'' clustering around those points which
gives rise to high-order accuracy. This change of variables is
analogous to that of the polar integration in~\citep{Kunyansky-1}, but
differs in the fact that it is applied on a rectangular mesh, hence
the ``rectangular-polar'' terminology we use to refer to this
methodology. The~\emph{sinh} transform~\citep{Johnston-1,Johnston-2}
was also tested as an alternative to the change of variables we
eventually selected: the latter method was preferred as the
\emph{sinh} change of variable does not appear to allow sufficient
control on the distribution of discretization points along the
integration mesh, which is needed in order to accurately resolve the
wavelength without use of excessive numbers of discretization points.

The proposed rectangular-polar approach, which yields high-order accuracy, leads
to several additional desirable properties, such as the use of Chebyshev
representations for the density, which allows the possibility to compute
differential geometry quantities needed for electromagnetic BIE by
means differentiation of Chebyshev series. Additionally, the nodes for Fej\'er's
first quadrature are the same as the nodes for the discrete orthogonality
property of Chebyshev polynomials, which make the computation of the Chebyshev
transforms straightforward. In addition to scattering by a bounded obstacle,
this integral equation solver can also be used in the context of the Windowed
Green function method for scattering by unbounded obstacles, as is the case of
layered media~\citep{Perez-1,Perez-2,Perez-3} and waveguide
problems~\citep{Garza-1}.

This paper is organized as follows. Section~\ref{Preliminaries} states
the mathematical problem of acoustic scattering and the integral
representations used in this paper. Section~\ref{Integration} then
describes the overall rectangular-polar integration strategy,
including details concerning the methodologies used to produce
integrals for smooth, singular and near-singular kernels as well as
edge-singular integral densities. A variety of numerical results for
open and closed scattering surfaces are then presented in
Section~\ref{numer}, emphasizing on the convergence properties of both
the forward map (which evaluates the action of the integral operator
for a given density) as well as the full scattering solver, and
demonstrating the accuracy, generality, and speed of the proposed
approach. Results of an application to a problem of scattering by a
geometry generated by CAD software is also presented in this section,
demonstrating the applicability to complex geometrical designs in
science and engineering. Section~\ref{Conclusions}, finally, presents
a few concluding remarks.

\section{Preliminaries\label{Preliminaries}} %------------------------------------------------
For conciseness, we consider the problem of acoustic scattering by a
sound-soft obstacle, though the methodology proposed is also
applicable to electromagnetic scattering and other integral-equation
problems involving singular kernels.

Let $\Omega$ denote the complement of an obstacle $D$ in
three-dimensional space, and let $\bdry$ denote the boundary of the
obstacle. Calling $\field$, $\field^\scat$ and $\field^\inc$ the
total, scattered and incident fields respectively, then
$\field = \field^\scat + \field^\inc$ obeys the Helmholtz equation
\begin{align}
  \Delta \field(\nex) + k^2 \field(\nex) = 0, \quad \nex \in \reals^3 \setminus \overline \bdry,
\end{align}
with wavenumber $k=2\pi/\lambda$, and the scattered field
$\field^\scat$ satisfies the Sommerfeld radiation condition and the
boundary condition
\begin{align}
  \field^\scat(\nex) = -\field^\inc(\nex), \quad \nex \in \bdry. 
\end{align}

As it is well known~\citep{Colton-1}, the scattered field can be represented in
terms of layer potentials, reducing the problem to a boundary integral equation
with singular kernels. The single- and double-layer potentials are defined by
\begin{align}
  \spot [\densi] (\nex) &= \int_\bdry G(\nex,\ney) \densi(\ney) \dsurf, \quad 
  \nex \in \reals^3 \setminus \overline \bdry, \\
  \dpot [\densi] (\nex) &= \int_\bdry \frac{\p G(\nex,\ney)}{\p \nn(\ney)} \densi(\ney) \dsurf, \quad 
  \nex \in \reals^3 \setminus \overline \bdry, 
\end{align}
respectively, and where
$G(\nex,\ney) = \exp{(i k |\nex-\ney|)} / (4\pi |\nex-\ney|)$ is the free-space
Green function of the Helmholtz equation, $\nn$ is the outwards-pointing normal
vector, and $\densi$ is the surface density. In this paper, we demonstrate the
proposed methodology by considering two problems under a unified scheme, namely
the case of closed and open surfaces.

\subsection{Closed surfaces}

For the case of a closed, bounded obstacle, we use a standard combined-field
formulation~\citep{Colton-1}
\begin{align}
  \field^\scat (\nex)  = \dpot [\densi] (\nex) - i k \spot [\densi] (\nex), \quad 
  \nex \in \reals^3 \setminus \overline \bdry,
\end{align}
which leads to the second-kind integral equation at the boundary
\begin{align}
  \displaystyle \frac{1}{2} \densi (\nex) + D[\densi](\nex) - 
  ik S[\densi](\nex) = - \field^\inc(\nex), \quad \nex \in \bdry,
\end{align}
where the single- and double-layer boundary operators are defined as
\begin{align}
  \displaystyle S[\densi](\nex) &= \int_\bdry G(\nex,\ney) \densi(\ney) \dsurf, \quad
  \nex \in \bdry, \\
  \displaystyle D[\densi](\nex) &= \int_\bdry \frac{\p G(\nex,\ney)}{\p \nn(\ney)} \densi(\ney) \dsurf, \quad
  \nex \in \bdry, 
\end{align}
respectively. 

This formulation is guaranteed to provide a unique density solution to the
scattering problem considered here~\citep{Colton-1}, and due to the nature of
this second-kind integral equation, the number of iterations for GMRES remains
essentially bounded as $k$ is increased.

\subsection{Open surfaces}

For scattering by an open surface, a combined field formulation is not possible
given that the jump conditions of the double-layer would imply different field
values from above and below the open surface. However, one can use a
single-layer formulation for such purpose, in this case we have  
\begin{align}
  \field^\scat (\nex) = \spot [\densi] (\nex), \quad 
  \nex \in \reals^3 \setminus \overline \bdry,
\end{align}
which leads to a first-kind integral equation
\begin{align}
  \displaystyle  
  S[\densi](\nex) = - \field^\inc(\nex), \quad \nex \in \bdry.
\end{align}
It is known that this formulation leads to an ill-conditioned system, which
manifests in increasing number of iterations for GMRES as $k$ becomes
larger. For the purpose of clarity, we choose to use this simple formulation to
present the open-surface scattering solver, which, as it will be shown in
Section~\ref{sec:disk}, was enough to achieve accuracies up to $10^{-13}$ in the
far field solutions for relatively low wavenumbers. For higher frequency
problems where the computational time may be a severe constrain, a regularized
version of this solver can be obtained by means of a composition of the
single-layer operator with the derivative of the double-layer,
see~\citep[Sec.~3]{Lintner-1}.

An important aspect of the open-surface case that leads to significant
difficulties is that the solution $\densi(\nex)$ presents a singularity of the
form
\begin{align}
  \densi \sim \frac{\Phi }{\sqrt{d}},
\end{align}
where $d$ is the distance to the edge and $\Phi$ is an infinitely differentiable
function throughout the boundary (including the
edge)~\citep{Lintner-1}. In~\citep{Lintner-1}, a strategy based on quadrature
rules for the exact singularity form where introduced, together with the polar
integration method from~\citep{Kunyansky-1}. We propose an alternative approach
in which, in addition of the polar-rectangular setup, a change of variables is
introduced in the parametrization of the surface. This change of variables is
such that its derivatives vanish at the edges, and thus smooths out the
integrands. Although not particularly designed to exactly match the singularity
of the density at edges, the proposed algorithm does provide a robust approach
for the treatment of the density-singularities that arise in the closed-surface
edge case (for which the degree of the singularity depends on the edge angle,
which may itself vary along the edge).

\subsection{Surface representation\label{surf}}

The proposed method assumes the scattering surface is described by a
set of $M$ non-overlapping ``logically-quadrilateral'' (LQ)
parametrized patches. This geometrical description is particularly
well suited for designs generated by CAD software, which generally can
export surface representations in terms of NURBS-based models---that
is, parametrizations expressed in terms of certain types of Rational
B-Splines. In fact, the potential afforded by direct use of
CAD-exported representations (without the expense, difficulty and
accuracy deterioration inherent in the use of surface triangulations)
provided the driving force leading to this paper: each NURBS trimmed
surface can be ``quadrilateralized'' without great difficulty, which
lends the method essentially complete geometric generality and a
remarkable ease of use.

In the proposed approach, then, the scattering surface $\bdry$ is
partitioned on the basis of a finite number $M$ of parametrizations
\[
\paramorig{q}: [-1,1]^2 \to \mathbb{R}^3 \quad (q = 1, 2,\dots, M),
\]
%\paramorig{q}(\uorig,\vorig) 
each one of which maps the unit square $[-1,1]^2$ in the
$(\uorig,\vorig)$-plane onto an LQ patch
within $\bdry$. Since we require the system of LQ patches to cover
$\bdry$, we have, in particular
\begin{align}
\displaystyle \bdry = \bigcup_{q=1}^M \left\{  \paramorig{q}(\uorig,\vorig) \; 
 \big| \; (\uorig,\vorig) \in [-1,1]^2 \right\}.
\end{align}

Clearly, any integral over $\bdry$ can be decomposed as a sum of
integrals over each of the patches. With the goal of treating the
closed- and open-surface cases in a unified way, we define the
boundary integrals
\begin{align}
  I(\nex) &= \sum_{q=1}^M I^{q}(\nex), \\
  I^{q}(\nex) &=  \int_{\bdry^q} \kernel(\nex,\ney) \densi(\ney) \dsurf, \label{eqn:iq}
\end{align}
where
\begin{align}
  \displaystyle \kernel(\nex,\ney) = 
  \begin{cases}
    \displaystyle \frac{\p G(\nex,\ney) }{\p \nn(\ney)} - ik G(\nex,\ney), & \text{(Closed surface)}, \\
    \displaystyle G(\nex,\ney), & \text{(Open surface)}.
  \end{cases}
\end{align}

In the following section, we propose a methodology for accurate
evaluation of the integrals $I^{q}(\nex)$ for a given density
$\densi(\ney)$. The solution to the integral equation problem is then
obtained via an application of the iterative linear-algebra solver
GMRES.

\section{Integration
  strategy\label{Integration}} %------------------------------------------------

The integration scheme we present consists of three main components:
(1)~Use of Fej\'er's first quadrature rule to compute integrals between
patches that are ``far'' away from each other, (2)~A rectangular-polar
high-order accurate quadrature rule for self-patch and near-patch
singular integrals, and (3)~A change of variables that resolves the
density singularities that arise at edges.

Using, for each $q$, the parametrization $\paramorig{q}$, the
integral~\eqref{eqn:iq} can be expressed in the form
\begin{align}
\label{eq:int1}
  I^q(\nex) &=  \int_{-1}^{1} \int_{-1}^{1} \kernelu^q(\nex,\uorig,\vorig) 
                \surfelu^q(\uorig,\vorig) \densiu^q(\uorig,\vorig) 
                \de \uorig \: \de \vorig,\quad (\nex \in \bdry),
\end{align}
where $\surfelu^q(\uorig,\vorig)$ denotes the surface Jacobian, and
where
\begin{align}
  \displaystyle \kernelu^q(\nex,\uorig,\vorig) &= 
  \kernel \left(\nex,\paramorig{q}(\uorig,\vorig) \right), \\
  \densiu^q(\uorig,\vorig) &= \densi \left(\paramorig{q}(\uorig,\vorig)\right).
\end{align}
The strategy proposed for evaluation of the integral in
equation~\eqref{eq:int1} depends on the proximity of the point $\nex$
to the $q$-th patch. For points $\nex$ that are ``far'' from the
patch, Fej\'er's first quadrature rule is used as detailed in
Section~\ref{sec:offpatch}. A special technique, the rectangular-polar
method, is then presented in Section~\ref{sec:close} to treat the case
in which $\nex$ is ``close'' or within the $q$-th patch. Before the
presentation of these smooth, singular and near-singular integration
methods, Section~\ref{sec:sing} describes the singular character of
integral-equation densities at edges, and proposes a methodology,
which is incorporated in the subsequent sections, for their treatment
in a high-order accurate fashion.

\subsection{Density singularities along edges}
\label{sec:sing}

The sharp edges encountered in general geometric designs have provided
a persistent source of difficulties to integral equation methods and
other scattering solvers.  The presence of such edges leads to
integrable singularities in the density solutions in both the
open-surface~\citep{Lintner-1} and
closed-surface~\citep{Costabel-1,Markkanen-1} cases. The strength of
the singularity, however, depends on the formulation and, for
closed-surfaces, on the angle at the edge, which is generally not
constant.

In order to tackle this difficulty in a general and robust manner, we
introduce a change of variables on the parametrization variables
$(\uorig,\vorig)$, a number of whose derivatives vanish along
edges. Such changes of variables can be devised from mappings such as
the one presented in~\citep[Sec. 3.5]{Colton-1}, which is given by
\begin{align}\label{colt_kress_1}
  \displaystyle \wck(\tau) = 2 \pi \frac{[v(\tau)]^p}{[v(\tau)]^p + [v(2 \pi-\tau)]^p}, 
  \quad  0\le \tau \le 2\pi,
\end{align}
where
\begin{align}\label{colt_kress_2}
  \displaystyle v(\tau) = \left( \frac{1}{p} - \frac{1}{2} \right) \left( \frac{\pi - \tau}{\pi}
  \right )^3 + \frac{1}{p} \left( \frac{\tau-\pi}{\pi} \right) + \frac{1}{2}.
\end{align}
It is easy to check that the derivatives of $\wck(\tau)$ up to order
$p-1$ vanish at the endpoints. The function $\wck(\tau)$ can then be
used to construct a change of variables to accurately resolve the edge
singularities while mapping the interval $[-1,1]$ to itself. The
change-of-variable mappings we use are given by
\begin{align} 
  \label{eq:change1}
  \uorig=\wu(\uu) =
  \begin{cases}
    \uu, & 
    \text{No edge on }\uorig \\
    -1 + \frac{1}{\pi} \wck\left(\pi[\uu+1] \right), & 
    \text{Edges at } \uorig \pm 1 \\
    -1 + \frac{2}{\pi} \wck \left(\frac{\pi}{2}[\uu+1] \right), & 
    \text{Edge at }\uorig=-1 \text{ only} \\
    -3 + \frac{2}{\pi} \wck \left(\pi + \frac{\pi}{2}[\uu+1] \right), & 
    \text{Edge at }\uorig=1 \text{ only}
  \end{cases}
\end{align}
and similarly
\begin{align} 
  \label{eq:change2}
  \vorig=\wv(\vv) =
  \begin{cases}
    \vv, & 
    \text{No edge on }\vorig \\
    -1 + \frac{1}{\pi} \wck \left(\pi[\vv+1] \right), & 
    \text{Edges at } \vorig \pm 1\\
    -1 + \frac{2}{\pi} \wck \left(\frac{\pi}{2}[\vv+1] \right), & 
    \text{Edge at }\vorig=-1 \text{ only} \\
    -3 + \frac{2}{\pi} \wck \left(\pi + \frac{\pi}{2}[\vv+1] \right), & 
    \text{Edge at }\vorig=1 \text{ only}
  \end{cases}
\end{align}
Incorporating the change of variables~\eqref{eq:change1}
and~\eqref{eq:change2}, the integral in equation~\eqref{eq:int1}
becomes an integral in which a weakly singular kernel is applied to a
finitely smooth function:
\begin{align}
  \label{eq:iq2}
  I^{q}(\nex) &= \int_{-1}^{1} \int_{-1}^{1} \kerneluu^q(\nex,\uu,\vv) 
                \surfeluu^q(\uu,\vv) \frac{d\wu}{d\uu}(\uu) \frac{d\wv}{d\vv}(\vv) 
                \densiuu^q(\uu,\vv) \de \uu \: \de \vv,\quad \nex \in \bdry,
\end{align}
where
\begin{align}
  \displaystyle \kerneluu^q(\nex,\uu,\vv) &= 
  \kernelu^q \left( \nex,\wu(\uu),\wv(\vv) \right), \\
  \densiuu^q (\uu, \vv) &= \densiu^q \left( \wu(\uu),\wv(\vv) \right), \\
  \param{q}(\uu, \vv) &= \paramorig{q} \left( \wu(\uu),\wv(\vv) \right), \\
  \surfeluu^q(\uu, \vv) &=\surfelu^q(\uorig,\vorig) .
\end{align}
(The edge-vanishing derivative factors in the integrand smooth-out any
possible edge singularities in the density
$\densiuu^q$~\cite{Colton-1}.) The proposed algorithm evaluates such
integrals by means of the ``smooth-density methods'' described in
Sections~\ref{sec:offpatch} and~\ref{sec:close} below.

\begin{figure}[H]
\centering
\includegraphics[scale=.4]{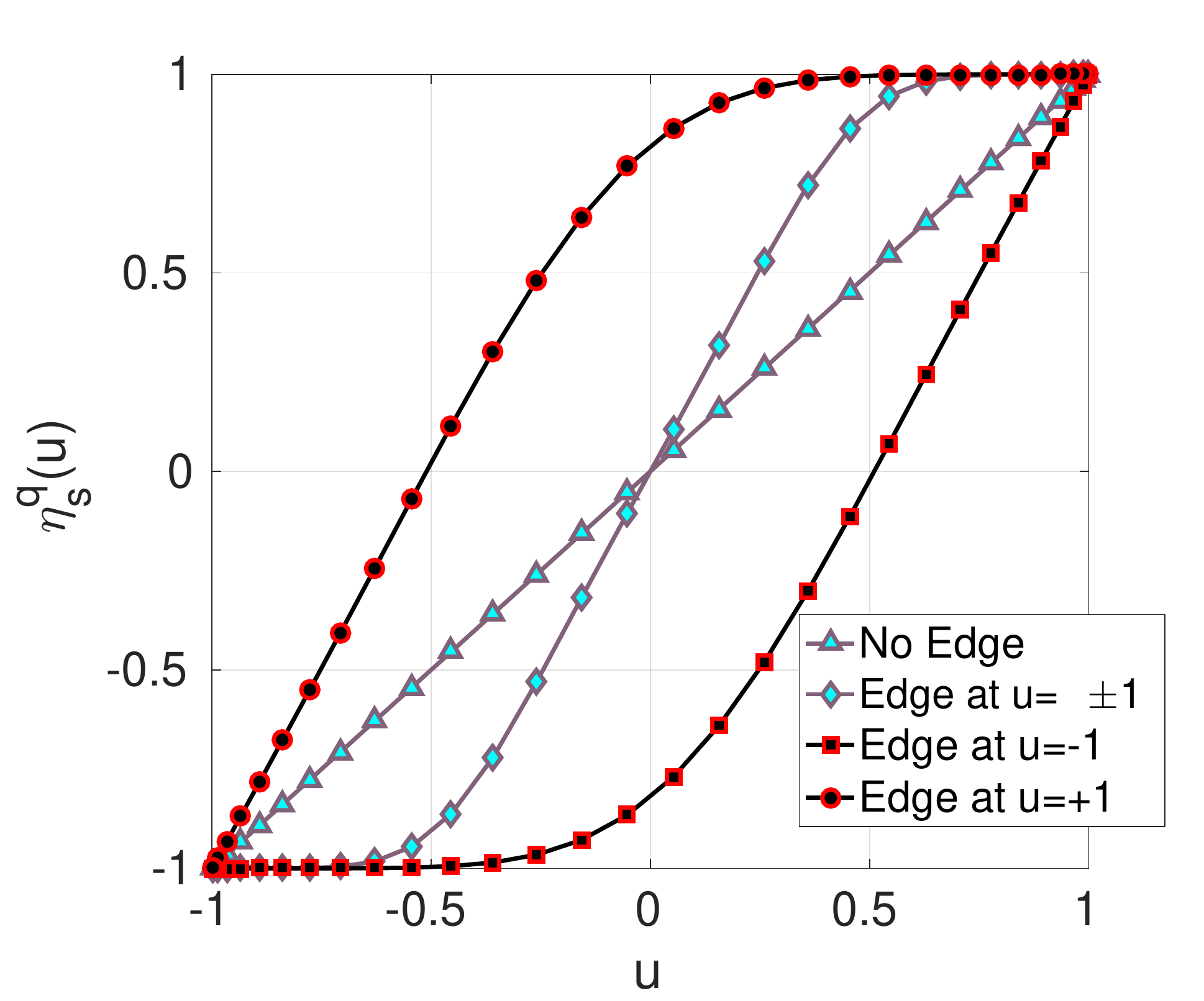}
\caption{Change of variables (equations~\eqref{eq:change1}
  and~\eqref{eq:change2}) used to resolve the density
  edge-singularities.\label{fig:edges}}
\end{figure}

\subsection{Non-adjacent integration}
\label{sec:offpatch}

The algorithm we use for the evaluation of the quantity $I^{q}(\nex)$,
defined by~\eqref{eq:int1}, is based on the
reformulation~\eqref{eq:iq2}---which takes into account all the
possible ways in which edge-singularities may or may not appear within an
integration patch. (The algorithm does assume that geometric
singularities may only appear along patch boundaries.)

In the ``non-adjacent'' integration case considered in this section,
in which the point $\nex$ is far from the integration patch, the
integrand in~\eqref{eq:iq2} is smooth---in view of the changes of
variables inherent in that equation, which, in particular, give rise
to edge-vanishing derivative factors that smooth-out any possible
edge-singularity in the density $\densiuu^q$ itself. (Using well known
asymptotics of edge singularities it is easy to check~\cite{Colton-1,BOT}
that the vanishing derivatives indeed smooth out all possible edge
singularities, to any desired order of smoothness, provided a
sufficiently high value of $p$ is used. Values of $p$ as low as $p=2$
are often found to be adequately useful, for accuracies of the order
of $1\%$, but values of $p$ of the orders of four to six or above can
enable significantly faster convergence and lower computing costs for
higher accuracies, as demonstrated e.g. in Figure~\ref{fig:cube}.)

In view of the smoothness of the integrands for the non-adjacent cases
considered presently ($\nex$ is far from the integration patch), then,
the integral in~\eqref{eq:iq2} could be accurately evaluated on the
basis of any given high-order quadrature rule. Our implementation
utilizes Fej\'er's first quadrature rule~\citep{Waldvogel-1}, which
effectively exploits the discrete orthogonality property satisfied by
the Chebyshev polynomials in the Chebyshev meshes used, and which also
allows for straightforward computation of the two-dimensional
Chebyshev transforms that are required as part of the singular and
near-singular integration algorithms described in
Section~\ref{sec:close}.

For a discretization using $N$ points, the quadrature nodes and
weights are given by
\begin{align}
  \displaystyle \fejernodes_j = \cos{ \left( \pi \frac{2j+1}{2N} \right) }, 
  \quad j=0, \dots, N-1,
\end{align}
\begin{align}
  \displaystyle \fejerw_j = \frac{2}{N} \left( 1 - 2 \sum_{\ell=1}^{\lfloor N/2
  \rfloor} \frac{1}{4\ell^2 - 1} \cos{( \ell \fejernodes_j )}  \right), 
  \quad j= 0, \dots, N-1,
\end{align}
respectively. Then using the Cartesian-product discretization
$\{ \uu_i=\fejernodes_i|i=0,\dots,\npts{\uu}{q}-1\}\times\{
\vv_j=\fejernodes_j|j=0,\dots,\npts{\vv}{q}-1\}$, the integral
in~\eqref{eq:int1} can be approximated by the quadrature expression
\begin{align}
  I^q(\nex) \approx \sum_{j=0}^{\npts{\vv}{q}-1} \sum_{i=0}^{\npts{\uu}{q}-1} 
  \kerneluu^q(\nex,\uu_i,\vv_j) 
  \surfeluu^q(\uu_i,\vv_j) \frac{d\wu}{d\uu} \frac{d\wv}{d\vv}\fejerw_i \fejerw_j \,\densiuu^q(\uu_i,\vv_j), 
  \quad \nex \in \setfar_q.
\end{align}
where $ \setfar_q$ represents the set of points that are
``sufficiently far'' from the $q$-th integration patch.

\subsection{Singular ``rectangular-polar'' integration algorithm and a
  new edge-resolved integral unknown}
\label{sec:close}
In this section we consider once again the problem of evaluation the
quantity $I^{q}(\nex)$, on the basis of the reformulation~\eqref{eq:iq2},
but this time (with reference to the last paragraph in
Section~\ref{sec:offpatch}) for points $\nex$ that are not
sufficiently far from the $q$-th integration patch. The set of all
such points, which contains only points that either lie within the
$q$-th patch or are ``close'' to it, will be denoted by $\setclose_q$.
The problem of evaluation of $I^{q}$ for $\nex \in \setclose_q$
presents a significant challenge in view of the singularity of the
kernel $\kernel(\nex,\ney)$ at $\nex = \ney$. 

In order to deal with this difficulty, we propose once again the use of
smoothing changes of variables: we will seek to incorporate in the
integration process a change of variables whose derivatives vanish at
the singularity or, for nearly singular problems, at the point in the
$q$-th patch that is closest to the singularity. In previous
implementations~\citep{Kunyansky-1,Lintner-1}, such changes of
variables required interpolation of the density $\densiuu^q$ from the
fixed nodes $(\uu_i,\vv_j)$ to the new integration points. The
interpolation step, though viable, can amount to a significant portion
of the overall cost. We thus propose,
instead, use of a precomputation scheme for which integrals of the
kernel times Chebyshev polynomials are evaluated with
high-accuracy. Since Chebyshev polynomials can easily be evaluated at
any point in their domain of definition, this approach does not
require an interpolation step. And, since these integrals are
independent of the density, they need only be computed once at the
beginning of any application of the algorithm, and reused in the
algorithm as part of any necessary integration processes in subsequent
linear-algebra (GMRES) iterations. Thus, for a given density
$\densiuu^q$ the overall quantity $I^{q}(\nex)$ with
$\nex\in\setclose_q$ can be computed by first obtaining the Chebyshev
expansion
\begin{align}
  \label{eq:series}
  \densiw^q(\uu,\vv) = \sum_{m=0}^{\npts{\vv}{q}-1} \sum_{n=0}^{\npts{\uu}{q}-1} 
  \densicheb_{n,m} T_n(\uu) T_m(\vv)  ,
\end{align}
of the modified edge-resolved (smooth) density
\[
  \densiw^q(\uu,\vv) = \frac{d\wu}{d\uu}(\uu) \frac{d\wv}{d\vv}(\vv)
  \densiuu^q(\uu,\vv),
\]
and then applying the precomputed integrals for Chebyshev densities.

In detail, the necessary Chebyshev coefficients $\densicheb_{n,m}$ can
be computed using the relation~\citep{num-recipies}
\begin{align} 
  \displaystyle \densicheb_{n,m} = \frac{\alpha_n \alpha_m}{\npts{\uu}{q} \npts{\vv}{q}}
  \sum_{j=0}^{\npts{\vv}{q}-1} \sum_{i=0}^{\npts{\uu}{q}-1}
  \densiw^q( \uu_i, \vv_j) T_n(\uu_i) T_m(\vv_j) 
\end{align}
that results from the discrete-orthogonality property enjoyed by
Chebyshev polynomials, where
\begin{align}
  \alpha_n = 
  \begin{cases}
    1, & n=0,\\
    2, & n\neq 0. 
  \end{cases}
\end{align}
As is well known, the Chebyshev coefficients $\densicheb_{n,m}$ can be
computed in a fast manner either by means of the FFT algorithm or, for
small expansion orders, by means of partial
summation~\citep[Sec. 10.2]{Boyd-1}. In practice, relatively small
orders and numbers of discretization points are used, and we thus
opted for the partial summation strategy.

Using the expansion~\eqref{eq:series} we then obtain
\begin{align}
  \displaystyle I^{q}(\nex) &= \int_{-1}^{1} \int_{-1}^{1} \kerneluu^q(\nex,\uu,\vv) 
                \surfeluu^q(\uu,\vv) \left( \sum_{m=0}^{\npts{\vv}{q}-1} \sum_{n=0}^{\npts{\uu}{q}-1} 
                \densicheb_{n,m} T_n(\uu) T_m(\vv) \right) \de \uu \: \de \vv
\end{align}
from which, exchanging the integrals with the sum, it follows that
\begin{align} 
  \label{eq:int2}
  \displaystyle I^{q}(\nex) = \sum_{m=0}^{\npts{\vv}{q}-1} \sum_{n=0}^{\npts{\uu}{q}-1} 
                \densicheb_{n,m} \int_{-1}^{1} \int_{-1}^{1} \kerneluu^q(\nex,\uu,\vv) 
                \surfeluu^q(\uu,\vv)   T_n(\uu) T_m(\vv) \, \de \uu \: \de \vv.
\end{align}
As mentioned above, the double integral is independent of the density,
and therefore it only depends on the geometry, the kernel, and the
target point $\nex \in \setclose_q$. For the computation of the
forward map, we only need to evaluate $I^{q}(\nex)$ for all {\em
  discretization points} $\nex \in \setclose_q$. Thus, in the proposed
strategy, the integral in~\eqref{eq:int2} must be precomputed for each
$q$ and for each combination of a target point $\nex \in \setclose_q$
and a relevant product of Chebyshev polynomials. Denoting the set of
all discretization points by
\begin{align}
  \displaystyle \setpts = \left\{ 
  \paramorig{q} \left( \wu(\uu_i), \wv(\vv_j) \right)
  \big| \; q = 1,\dots,M, \; i=1,\dots,\npts{\uu}{q}, \; j=1,\dots,\npts{\vv}{q}
 \right  \},
\end{align}
and using the weights
\begin{align}
  \label{eq:int3}
  \displaystyle \chebyw{q}{n}{m}{\ell} = \int_{-1}^{1} \int_{-1}^{1} 
  \kerneluu^q(\nex_\ell,\uu,\vv) 
  \surfeluu^q(\uu,\vv)   T_n(\uu) T_m(\vv) \, \de \uu \: \de \vv, \quad
  \text{for each } \nex_\ell \in \left\{ \setpts \cap \setclose_q \right\},
\end{align}
equation~\eqref{eq:int2} becomes
\begin{align} 
  \displaystyle I^{q}(\nex_\ell) = \sum_{m=0}^{\npts{\vv}{q}-1} \sum_{n=0}^{\npts{\uu}{q}-1} 
  \densicheb_{n,m} \, \chebyw{q}{n}{m}{\ell} \;.
\end{align}

We now turn our attention to the accurate evaluation of the integrals
in equation~\eqref{eq:int3}. The previous method~\citep{Kunyansky-1}
utilizes (in a different context, and without precomputations) a polar
change of variables that cancels the kernel singularity and thus gives
rise to high-order integration.  Reference~\citep{Kunyansky-1} relies
on overlapping parametrized patches and partitions of unity to
facilitate the use of polar integration.  In the case in which
non-overlapping LQ patches are utilized, the use of polar integration
requires design of complex quadratures near all patch
boundaries~\cite{Lintner-1}. To avoid these difficulties, we propose
use of certain ``rectangular-polar'' changes of variables which, like
the edge changes-of-variables utilized in Section~\ref{sec:offpatch},
is based on use of the
functions~\eqref{colt_kress_1}--\eqref{colt_kress_2} for suitable
values of $p$.

We thus seek to devise rectangular-polar change of variables that can
accurately integrate the kernel singularity, either in the self-patch
problem, for which the singularity lies on the integration patch and
for which changes of variables should have vanishing derivatives at
the target point $\nex_\ell$, or in the near-singular case, for which
vanishing change-of-variable derivatives should occur at the point in
the $q$-th patch that is closest to the observation point $\nex$. To
achieve this, it is necessary to consider the value
\begin{align}\label{min_proj}
  \displaystyle \left( \uproj{q}{\ell},\vproj{q}{\ell} \right) = 
  \operatorname*{arg\,min}_{(u,v)\in[-1,1]^2}
  \left| \, \nex_\ell - \param{q}(\uu,\vv) \, \right|, 
\end{align}
which can be found by means of an appropriate minimization
algorithm. In view of its robustness and simplicity, our method
utilizes the golden section search algorithm~(see
\citep[Sec. 10.2]{num-recipies}) for this purpose, with initial bounds
obtained from a direct minimization over all of the the original
discretization points $\nex_\ell$ in the patch.  Relying on the
coordinates~\eqref{min_proj} of the projection point in the
near-singular case, and using the same notation
$\left( \uproj{q}{\ell},\vproj{q}{\ell} \right)$ for the coordinates
of the singular point in the self-patch problem, the relevant
rectangular-polar change-of-variable can be constructed on the basis
of the one-dimensional change of variables
\begin{align}
  \label{eq:sing-cv}
  \displaystyle \wsingu{\alpha}(\tau) = 
  \begin{cases}
    \displaystyle \alpha + \left( \frac{\text{sgn}{(\tau)} -
        \alpha}{\pi} \right)
    \wck \left( \pi |\tau| \right), & \text{for } \alpha \not= \pm 1,  \\
    \displaystyle \alpha - \left( \frac{ 1 + \alpha}{\pi} \right)
    \wck \left( \pi \left| \frac{\tau-1}{2} \right| \right), & \text{for } \alpha =  1,  \\
    \displaystyle \alpha + \left( \frac{ 1 - \alpha}{\pi} \right) \wck
    \left( \pi \left| \frac{\tau+1}{2} \right| \right), & \text{for }
    \alpha = -1.
  \end{cases}
\end{align}

Indeed, a new use of Fej\'er's first quadrature rule now yields
\begin{align}
  \label{eq:singcheby}
  \displaystyle \chebyw{q}{n}{m}{\ell} \approx \sum_{j=1}^{\nsing^v}
  \sum_{i=1}^{\nsing^u} 
  \kerneluu^q \left( \nex_\ell, \uu_i^{q,\ell}, \vv_j^{q,\ell} \right) 
  \surfeluu^q \left( \uu_i^{q,\ell}, \vv_j^{q,\ell} \right)
  T_n \left( \uu_i^{q,\ell} \right)
  T_m \left( \vv_j^{q,\ell} \right) \,
  \wsinguweight_i^{\uu,q,\ell} \wsinguweight_j^{\vv,q,\ell}
  \fejerw_i \, \fejerw_j
\end{align}
where 
\begin{align}
  \displaystyle \uu_i^{q,\ell} &= \wsingu{ \uproj{q}{\ell} } 
  \left( \fejernodes_i \right), \quad \text{for } i = 1,\dots, \nsing^u, \\
  \displaystyle \vv_j^{q,\ell} &= \wsingu{ \vproj{q}{\ell} } 
  \left( \fejernodes_j \right), \quad \text{for } j = 1,\dots, \nsing^v,
\end{align}
are the new quadrature points, and where
\begin{align}
  \displaystyle \wsinguweight_i^{\uu,q,\ell} &= \frac{d\wsingu{\uproj{q}{\ell}}}{d\tau}
  \left( \fejernodes_i \right), \quad \text{for } i = 1,\dots, \nsing^u, \\
  \displaystyle \wsinguweight_j^{\vv,q,\ell} &= \frac{d\wsingu{ \vproj{q}{\ell}}}{d\tau} 
  \left( \fejernodes_j \right), \quad \text{for } j = 1,\dots, \nsing^v,
\end{align}
denote the corresponding change-of-variable weights.  Using sufficiently large
numbers $\nsing^u$ and $\nsing^v$ of discretization points along the $u$ and $v$
directions to accurately resolve the challenging integrands, all singular and
near-singular problems can be treated with high accuracy under discretizations
that are not excessively fine (see Figure~\ref{fig:conv}).

\begin{figure}[H]
\centering
\includegraphics[scale=.4]{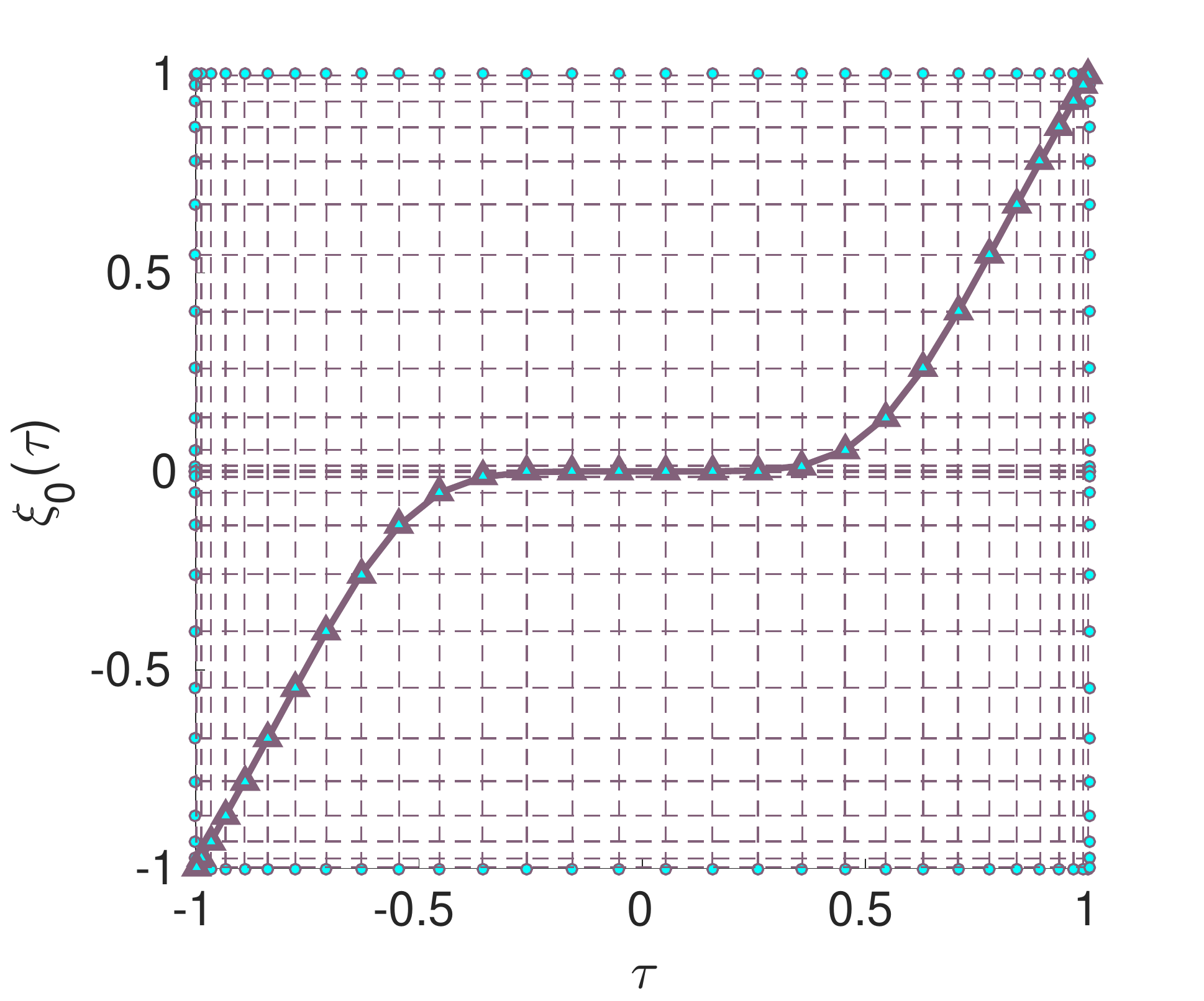}
\includegraphics[scale=.4]{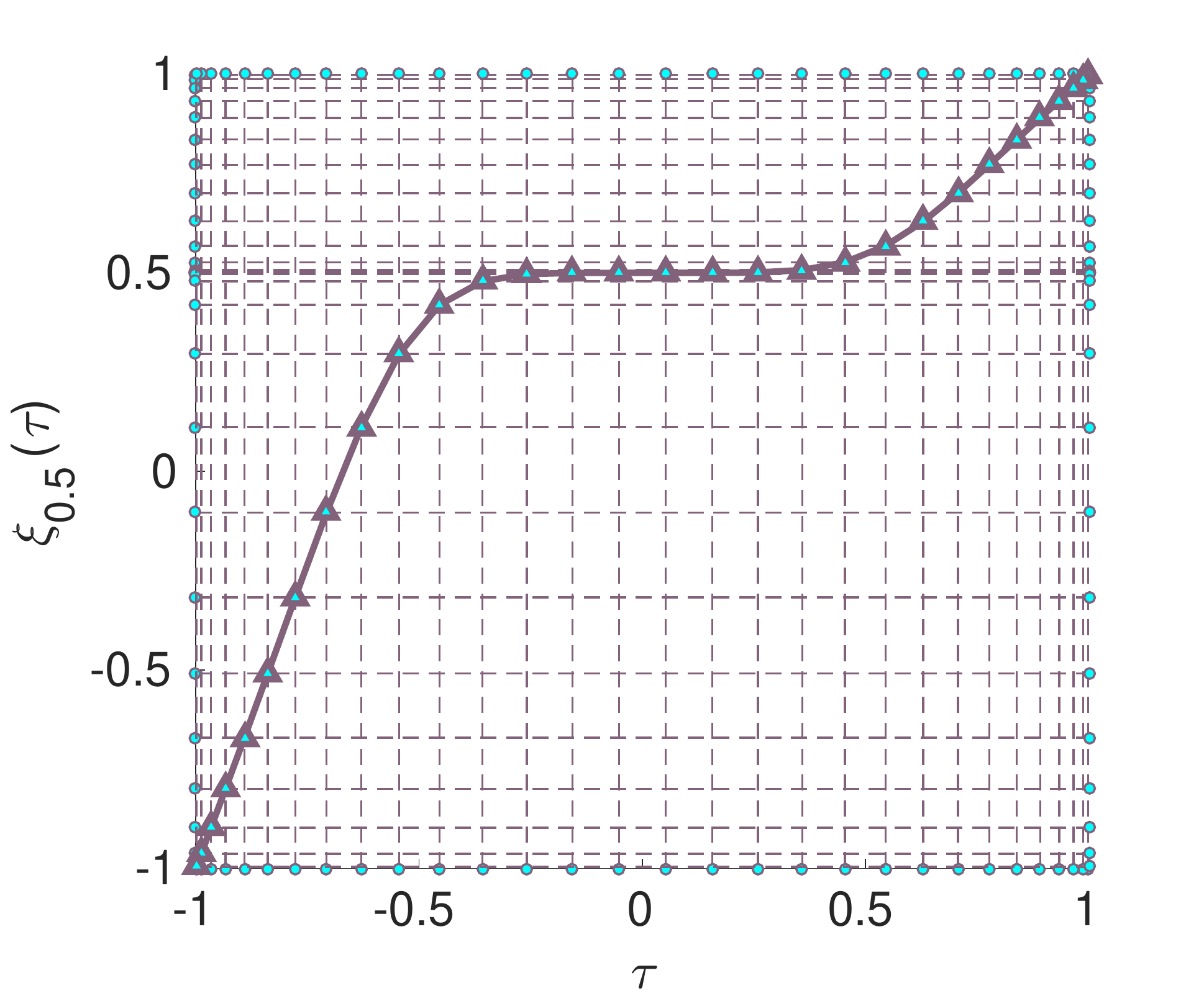} \\
\includegraphics[height=2.5in]{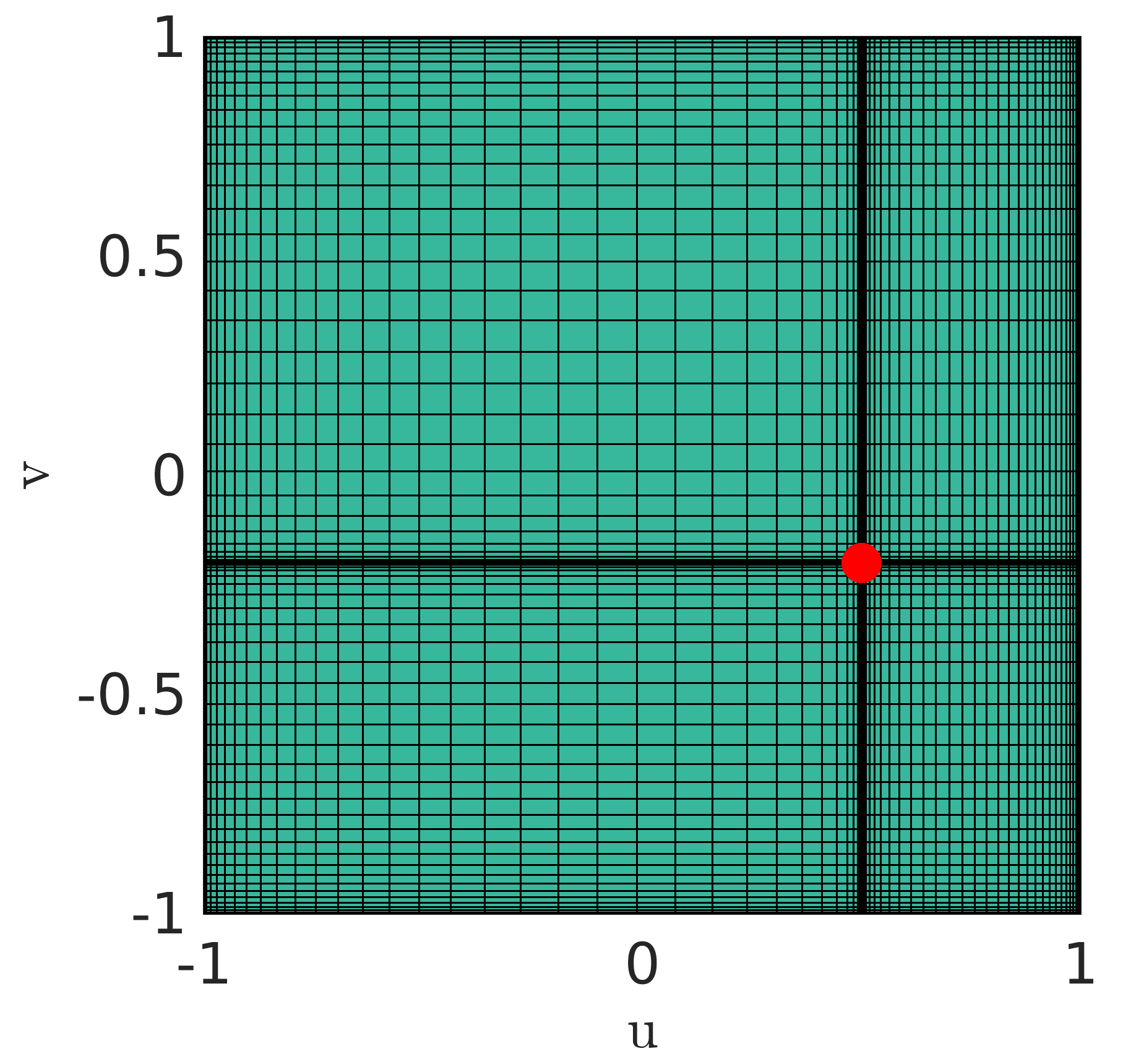} \hspace{0.5in}
\includegraphics[height=2.5in]{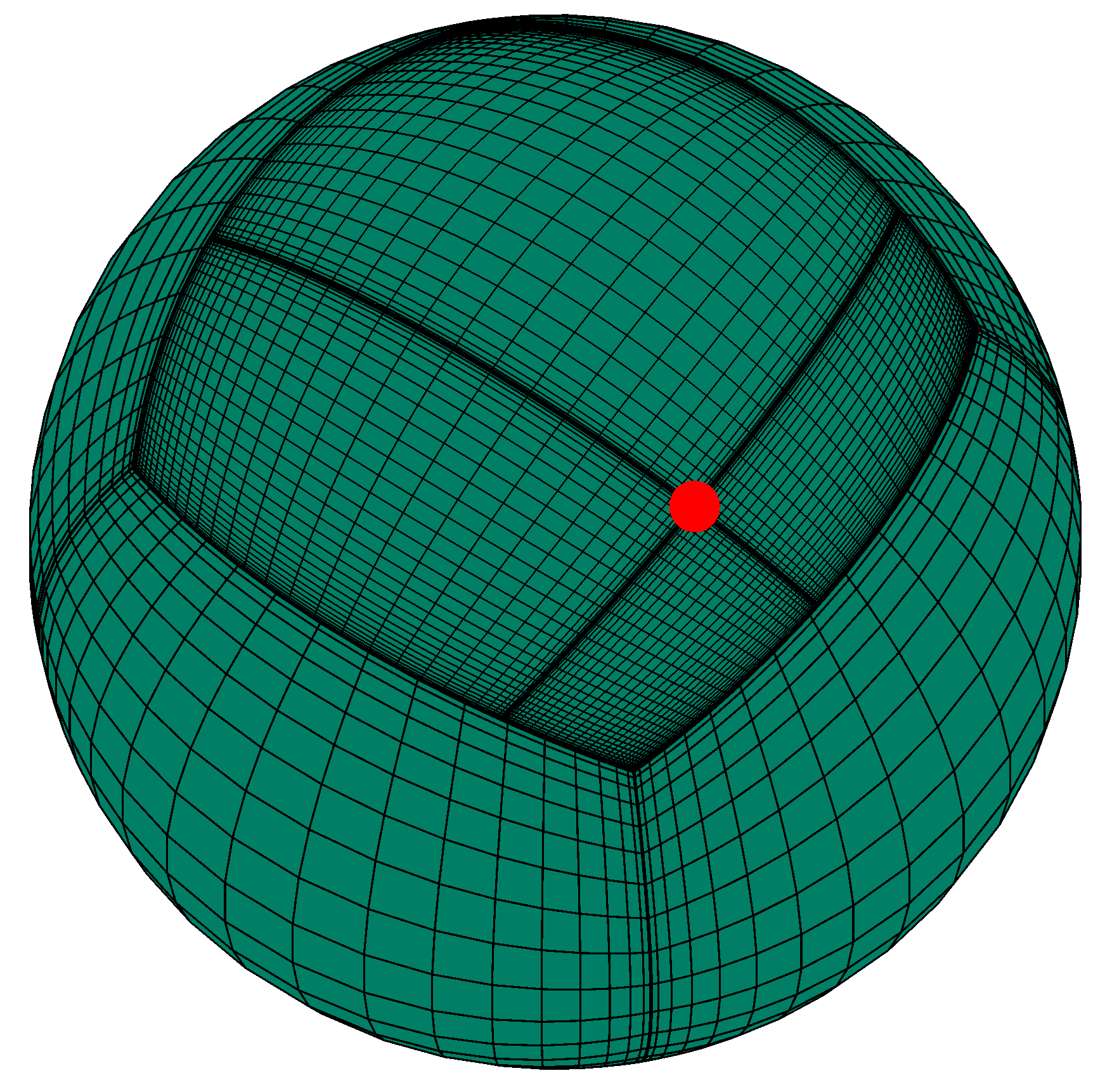}
\caption{Top: Changes of variables (equation~\eqref{eq:sing-cv}) used
  to resolve the kernel singularity for two different values of
  $\alpha$. Bottom: Mesh (in both parameter space [left] and real
  space [right]) produced by the rectangular-polar change of variables
  to resolve the kernel singularity located at the point marked in
  red. \label{fig:cv_sing}}
\end{figure}

\subsection{Computational cost}
Let us now estimate the computational cost for the
proposed method, focusing on the singular and near-singular integrals,
as the cost of the non-adjacent interactions arises trivially from a
double sum, and can be accelerated by means of either an equivalent
source scheme~\citep{Kunyansky-1,Kunyansky-2} or by a fast multipole
(FMM) approach~\citep{Gumerov-1}.

For the purposes of our computing-time estimates, let $N$ denote the maximum of
the one dimensional discretization sizes $\npts{\uu}{q}$ and $\npts{\vv}{q}$
over all patches ($1\leq q\leq M$), and let $N'_\text{close}$ denote the
maximum, over all the patches, of the numbers of points that are close to the
patch (they lie in $\setclose_q$), but which are not contained in the $q$-th
patch. Additionally, let $\nsing^u=\nsing^v=\nsing$ denote the number of
quadrature points used for singular precomputations. With these notations we
obtain the following estimates in terms of the bounded integer values $N$ (in
the range of one to a few tens); the (large, proportional to the square of the
frequency, for large frequencies) number $M$ of patches, and the related
bounded parameters $\nsing$ (of the order of one to a few hundreds):
\begin{itemize}
\item Cost of precomputations:
  $\bigo( M \nsing^2 ( N^2 + N'_\text{close} ) ) $ operations.
  \item Cost of forward map:
  \begin{itemize}
    \item Chebyshev transform (partial summation): $\bigo( M N^3 )$.
    \item Singular and near-singular interactions: $\bigo( M N^2 ( N^2 + N'_\text{close} ) )$.
    \item Non-adjacent interactions $\bigo( (M-1)^2 N^4 )$ (or
      $M^\alpha N^4$ with $\alpha$ significantly smaller than two if
      adequate acceleration algorithms are utilized).
  \end{itemize}
\end{itemize}

\subsection{Patch splitting for large problems}
\label{sec:mem}
Each patch requires creation and storage of a set of self-interaction weights
$\chebyw{q}{n}{m}{\ell}$, for $q=1,\dots,M$, $n=1,\dots,N$, $m=1,\dots,N$ and
$\ell=1,\dots,N^2$, at a total storage cost of $\bigo(M N^4)$ double-precision
complex-valued numbers. Additionally, weights also need to be stored for the
$N'_\text{close}$ near-singular points for each patch, and are dependent on the
target point, then the total storage for the singular and near-singular weights
is $\bigo ( M N^4+ M N^2 N'_\text{close})$.

In order to eliminate the need to evaluate and store a large number of
weights that result as $N$ is increased, it is possible to instead
increase the number of patches $M$---which causes the necessary number
of weights to grow only linearly. In these regards it is useful to
consider the following rule of thumb: in practice, as soon as the
wavelength is accurately resolved by the single-patch algorithm, due
to the spectral accuracy of Fej\'er's first quadrature, only a few
additional points per patch are needed to produce accuracies of the
order of several digits. In view of the estimates in this and the
previous section, parameter selections can easily be made by seeking
to optimize the overall computing time given the desired accuracy and
available memory.

\section{Numerical results\label{numer}} %--------------------------------------

This section presents a variety of numerical examples demonstrating
the effectiveness of the proposed methodology. The particular
implementation for the numerical experiments was programmed in Fortran
and parallelized using OpenMP. The runs were performed on a single
node of a dual socket Dell R420 with two Intel Xenon E5-2670 v3 2.3
GHz, 128GB of RAM. Unless otherwise stated, all runs where performed
using 24 cores.

\subsection{Forward map convergence}
The accuracy of the overall solver depends crucially on the accuracy
of the forward map computation. In this section we verify that the
proposed methodology yields uniformly accurate evaluations of the
action of the integral operator throughout the surface of the
scatterer. To do so, we consider the eigenfunctions and eigenvalues of
the single- and double-layer operators for Helmholtz
equation~\citep[Sec. 3.2.3]{Nedelec-1}:
\begin{align}
  \soper [Y_\ell^m(\theta,\varphi)] &= k j_\ell(k) h_\ell^{(1)}(k) 
  Y_\ell^m(\theta,\varphi), \\
  \displaystyle \doper [Y_\ell^m(\theta,\varphi)] &= \frac{k^2}{2} 
  \left[ j_\ell(k) \frac{\de}{\de k} h_\ell^{(1)}(k) + h_\ell^{(1)} 
  \frac{\de}{\de k} j_\ell(k) \right] 
  Y_\ell^m(\theta,\varphi),
\end{align}
where $j_\ell(k)$ and $h_\ell^{(1)}(k)$ are the spherical Bessel
function of the first kind and spherical Hankel function respectively,
and where $Y_\ell^m(\theta,\varphi)$ are the spherical harmonics. (For
the spherical Hankel function $h_\ell^{(1)}(z)$ we have used the
convention in~\citep{Nedelec-1}:
$h_\ell^{(1)}(z) = - y_\ell(z) + i \; j_\ell(z)$, where $y_\ell$ is
the $\ell$-th Neumann function.) 

\begin{figure}[H]
\centering
\includegraphics[height=2.6in]{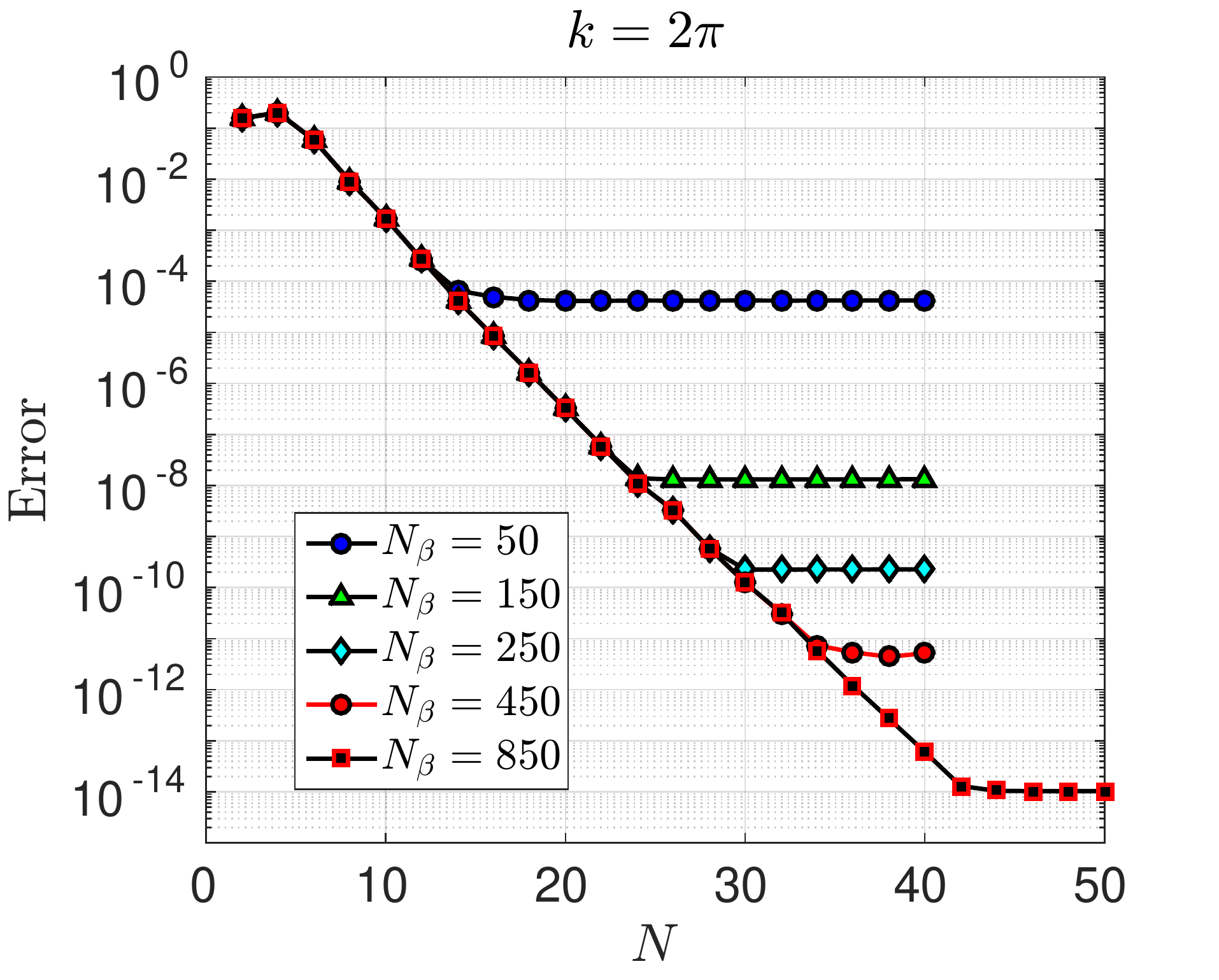}
\includegraphics[height=2.6in]{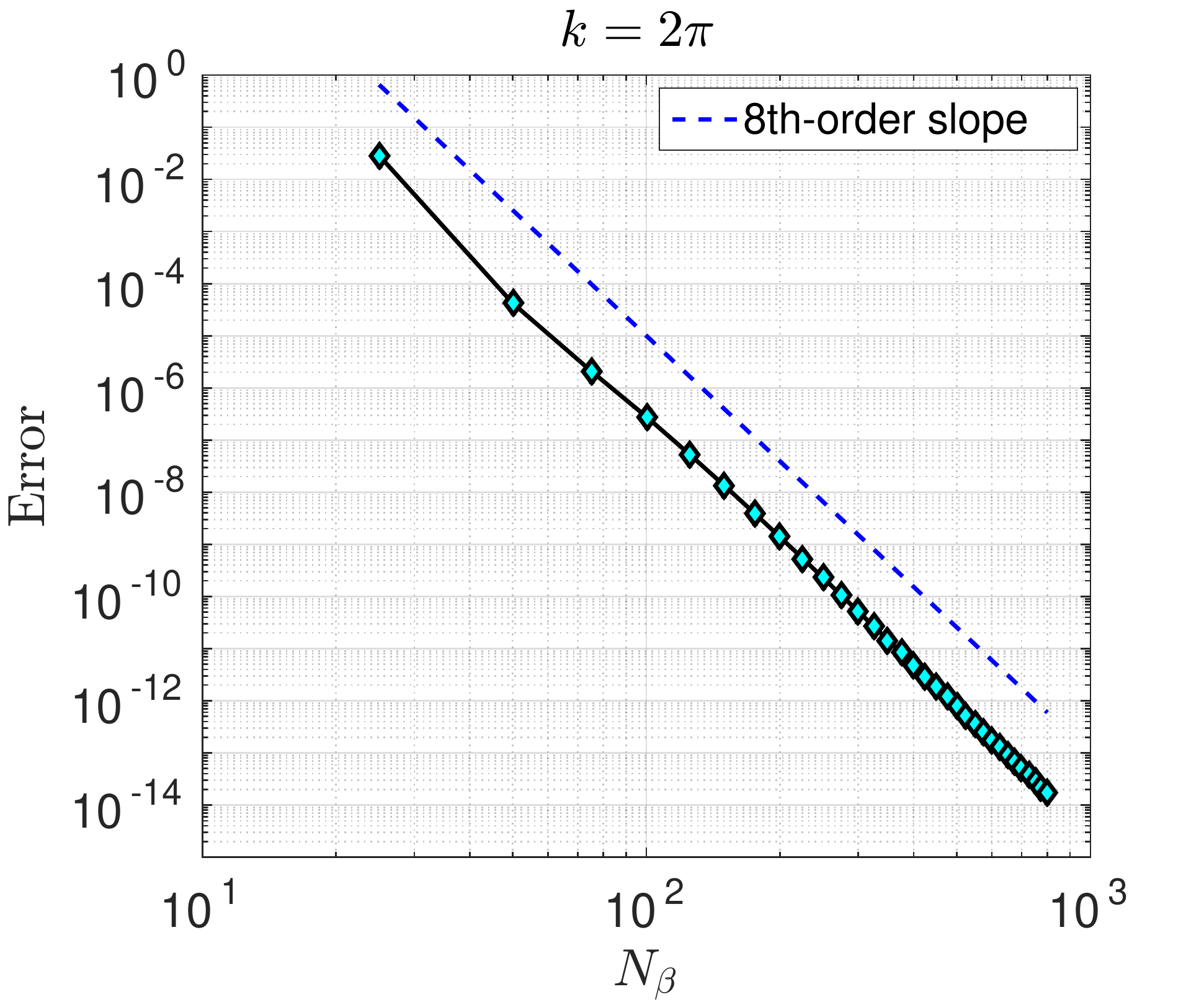}\\
\includegraphics[height=3in]{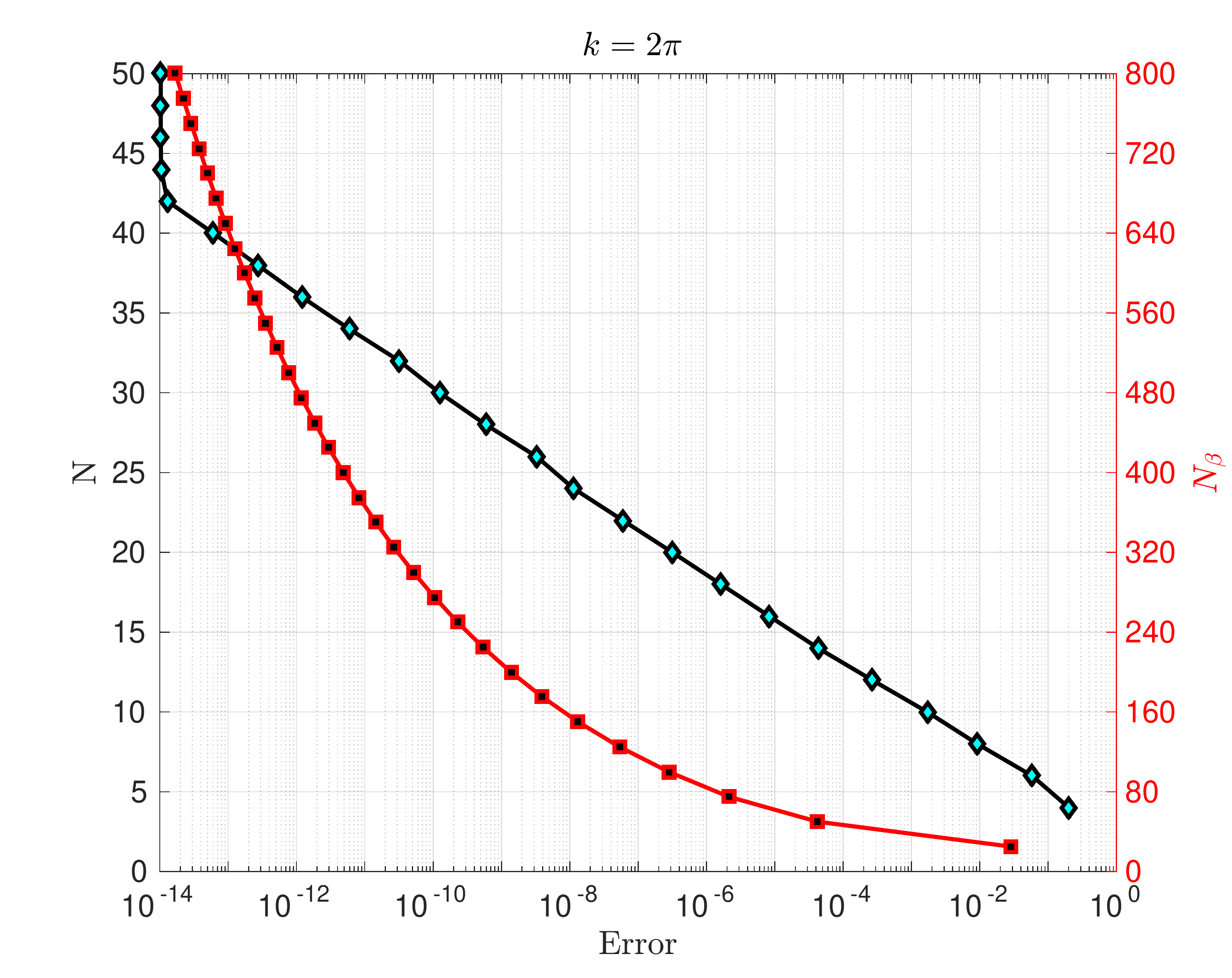}
\caption{Convergence on the forward map of the combined field
  formulation for a unit sphere and wavenumber $k=2\pi$. Top-left:
  Error as we increase the number of points per patch (per dimension)
  for different values of $N_\beta$ used in the
  precomputations. Top-right: Error as a function of
  $N_\beta$. Bottom: Values of $N$ (black curve) and $N_\beta$ (red
  curve) for a given error.\label{fig:conv}}
\end{figure}

\begin{figure}[H]
\centering
\includegraphics[scale=0.4]{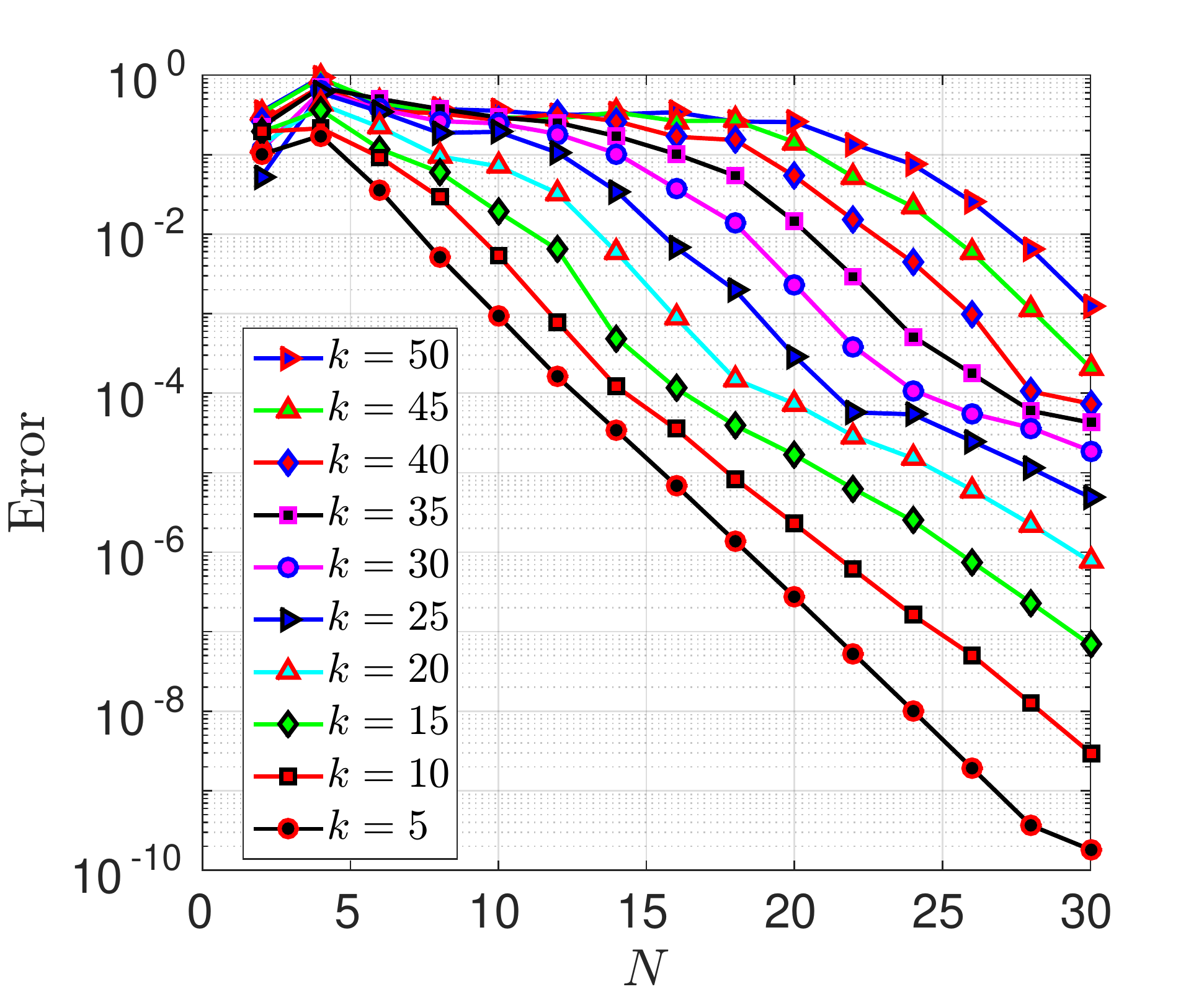} 
\caption{Error in the forward map of the combined field formulation for a unit
  sphere and for different wavenumbers as we increase the number of points per
  patch while keeping the numbers of patches fixed. A value of $\nsing=250$ was
  used to guarantee that the error in the singular precomputations is not larger
  than the wavelength discretization error.\label{fig:conv_k}}
\end{figure}

\begin{figure}[H]
\centering
\includegraphics[width=2in]{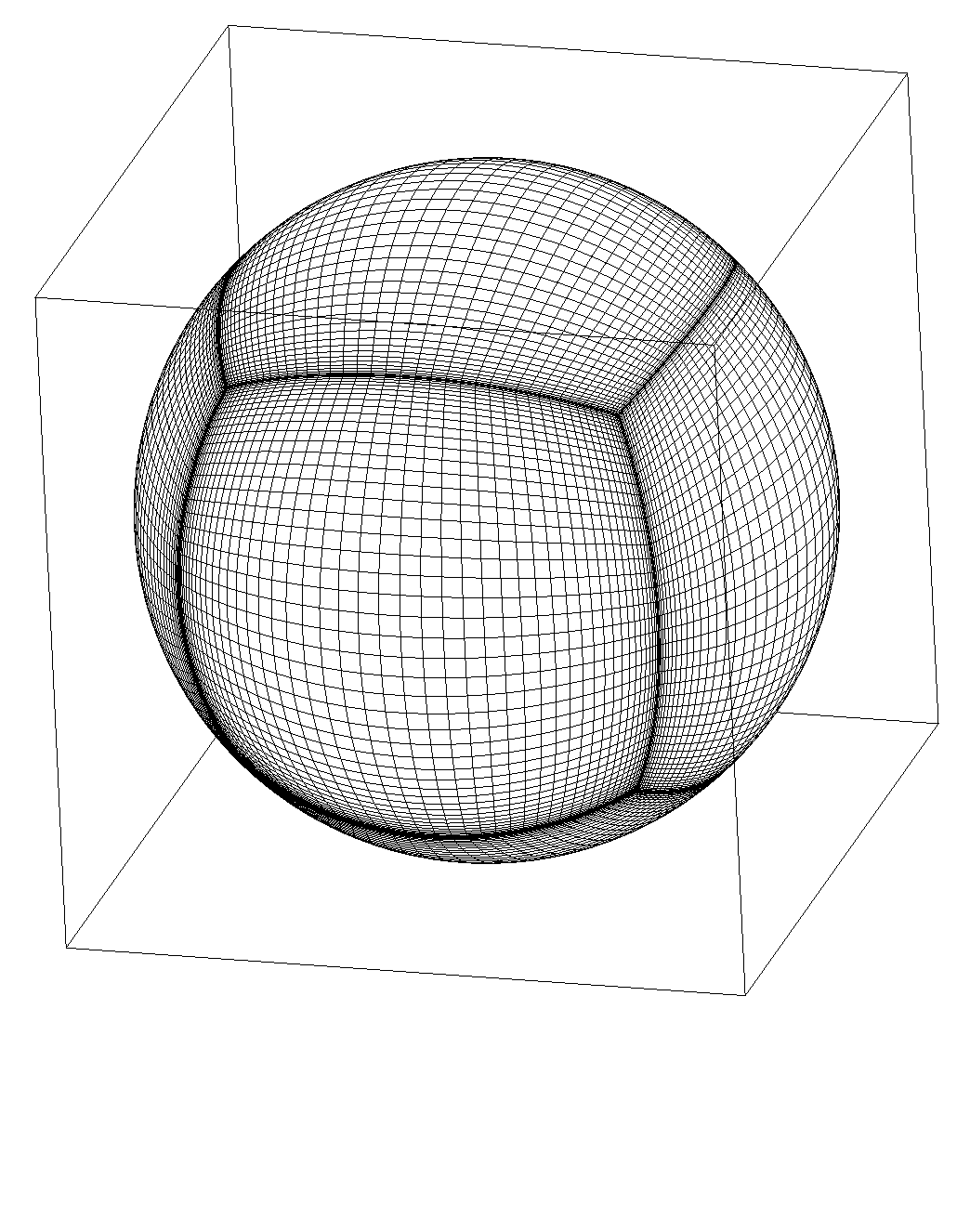} 
\includegraphics[width=2in]{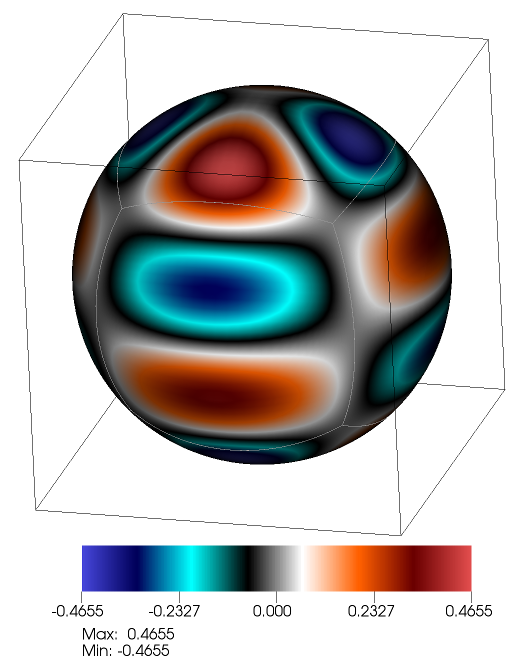} 
\includegraphics[width=2in]{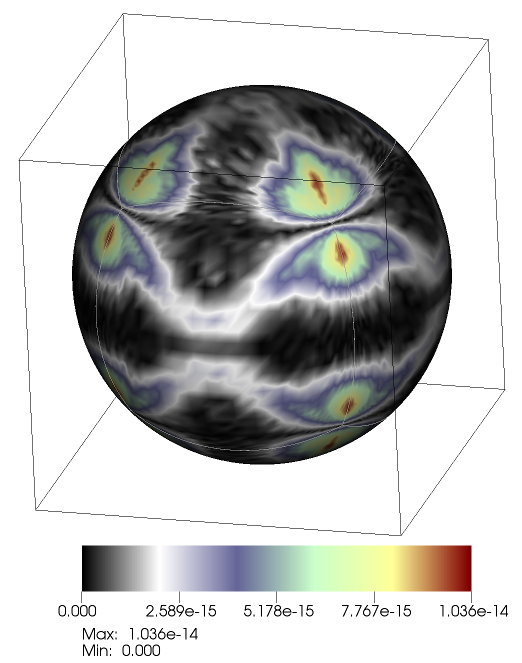}
\caption{Mesh, forward map and pointwise error (from left to right)
  for a spherical harmonic density $(5,2)$ using a mesh of $N=50$ and
  $N_\beta=850$; the error is essentially uniform and accuracies
  close to machine precision are achieved.\label{fig:error_fwd}}
\end{figure}

Figures~\ref{fig:conv} and~\ref{fig:conv_k} present convergence
results for the combined field formulation we use, and in particular
we demonstrate that the method is capable of obtaining accuracies
close to machine precision (Figure~\ref{fig:error_fwd}).

\begin{figure}[H]
\centering
\includegraphics[width=2in]{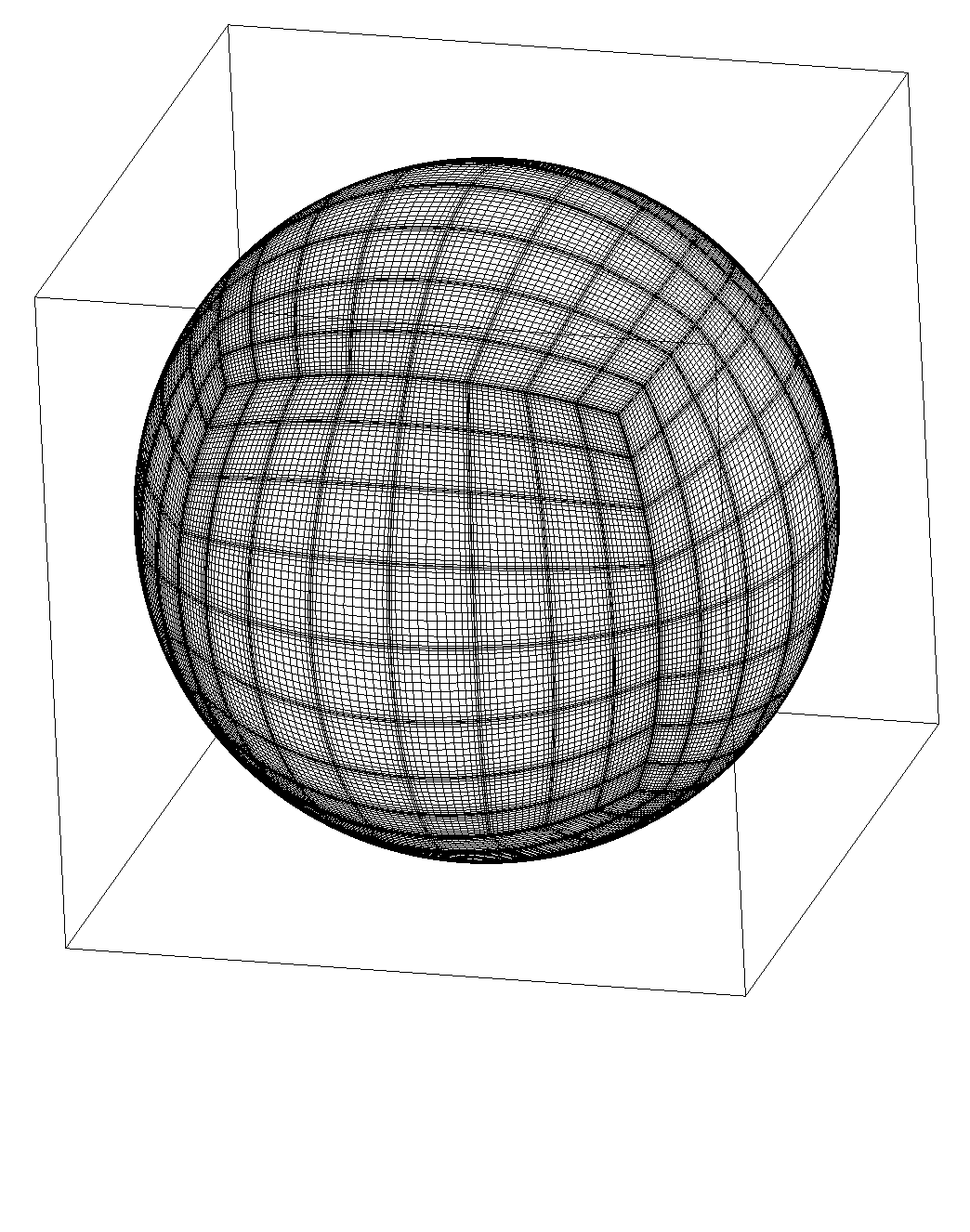} 
\includegraphics[width=2in]{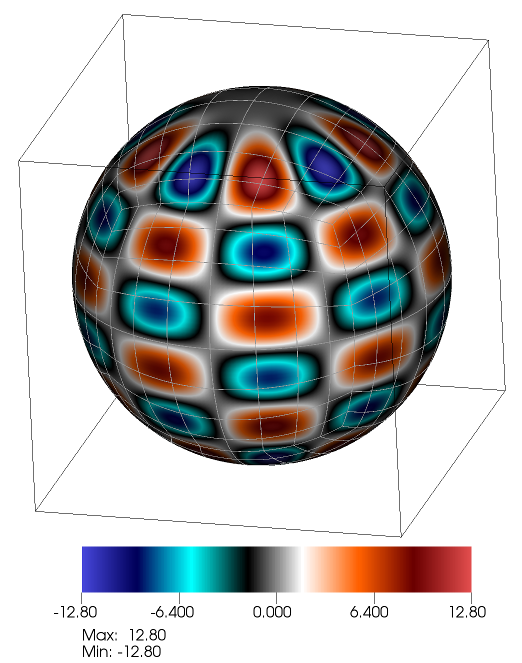} 
\includegraphics[width=2in]{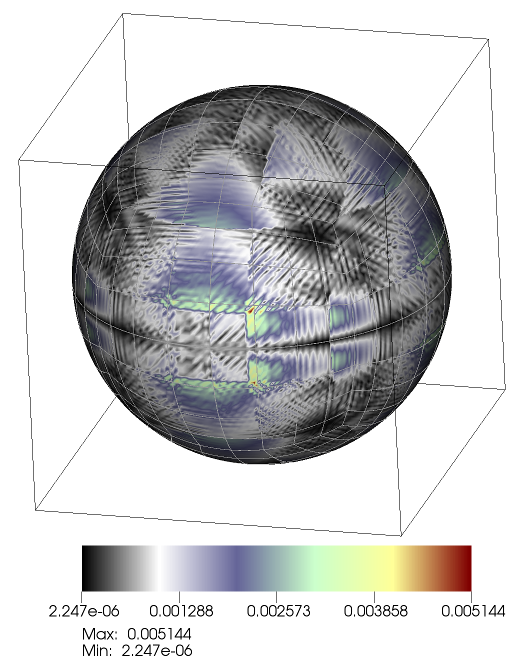}
\caption{Mesh, forward map and pointwise error (from left to right)
  for a spherical harmonic density $(10,5)$. The figures show the
  patch splitting strategy for high frequencies. In this case $k=100$
  which corresponds to scattering by a sphere $31.8\lambda$ in
  diameter, and and $6\times 8 \times 8$ patches ($8\times 8$
  subpatches in each one of $6$ initial patches) with
  $N=14$.\label{fig:error_sub}}
\end{figure}

As indicated in Section~\ref{sec:mem}, for high-frequency problems it is
beneficial to split the patches into smaller ones rather than increasing the
numbers of points per patch, given that the storage only grows linearly as the
number of patches is increased while keeping the number of points per patch
constant. In order to determine the optimal balance between accuracy and
efficiency, it must be considered that there are two factors that determine the
accuracy of the method: (1)~The order $N$ of the Chebyshev expansions used
(i.e. the number of points per patch per dimension), and (2)~The number of
points per wavelength. Figure~\ref{fig:error_sub} displays the pointwise error
in the forward map for a high frequency case, and Table~\ref{tab:tab1} presents
test results for several simple patch-splitting configurations, where the number
of points per wavelength is calculated by the formula
\begin{align}
  \text{Points per } \lambda = \frac{N}{L/\lambda},
\end{align}
where $L^2={4\pi/M}$ is the average area of the quadrilateral patches for the
sphere.

\begin{table}
  \centering
  \begin{tabular}{| c | c | c | c | r | c | c | c | }
    \hline
    $N$  & $\nsing$ & Patches & Points per $\lambda$ & Unknowns & Time (prec.) & Time (1 iter.) & Error \\ \hline
    $8$  & $50$ & $6 \times 5 \times 5 $   & $1.7$ & $9600$   & $1.39$ s  & $0.18$ s & $71.3\%$  \\ \hline
    $12$ & $60$ & $6 \times 5 \times 5 $   & $2.6$ & $21600$  & $3.40$ s  & $0.83$ s & $2.16\%$  \\ \hline
    $16$ & $80$ & $6 \times 5 \times 5 $   & $3.5$ & $38400$  & $9.26$ s  & $2.41$ s & $0.0814\%$  \\ \hline
    $8$  & $50$ & $6 \times 10 \times 10 $ & $3.5$ & $38400$  & $16.74$ s  & $2.78$ s  & $0.336\%$  \\ \hline
    $12$ & $70$ & $6 \times 10 \times 10 $ & $5.2$ & $86400$  & $47.26$ s  & $13.01$ s & $0.0238\%$  \\ \hline
    $16$ & $90$ & $6 \times 10 \times 10 $ & $6.9$ & $153600$ & $126.09$ s & $40.29$ s & $0.000355\%$  \\ \hline
  \end{tabular}
  \caption{Errors in the forward map (relative to the maximum forward map value)
    of the combined field operator for various patch splitting configurations
    and a spherical harmonic density (5,2).  For the results in this table,
    $k = 100$, a sphere of diameter $31.8\lambda$ was used, and all times reported
    for the precomputations and forward map where obtained using 24
    computing cores. At this frequency, $N=16$ suffices to produce an
    accuracy of $0.08\%$ on the basis of six $5\times 5$
    patches.\label{tab:tab1}}
\end{table}

\subsection{Edge geometries}

\begin{figure}[H]
\centering
\includegraphics[scale=0.4]{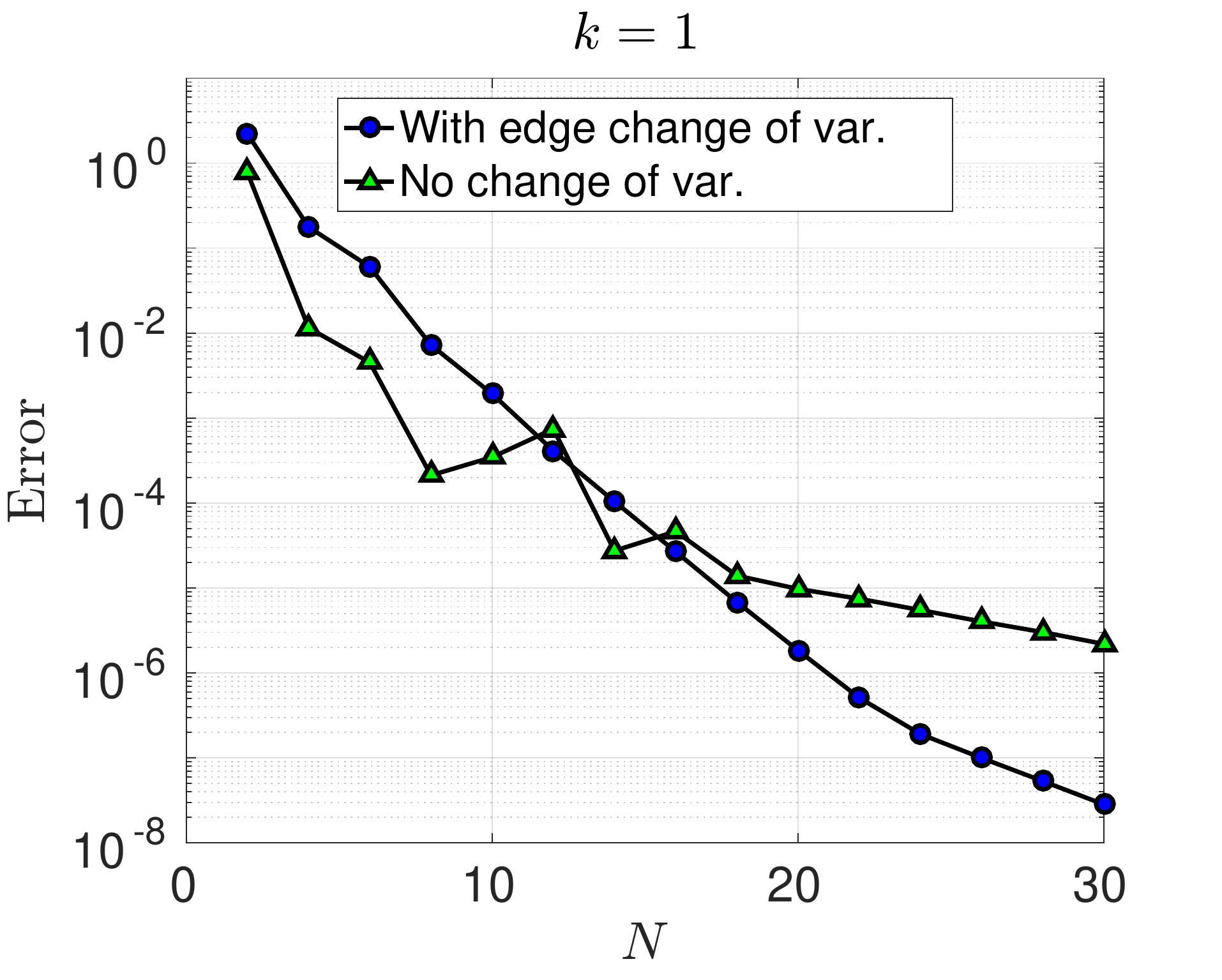}
\caption{Maximum (absolute) far-field error for the problem of scattering by a
  cube of size $2\times 2 \times 2$ and with $k=1$. The plot shows results
  excluding use of an edge change of variables (in green) and including an edge
  the change of variables with $p=2$ (in blue). The maximum value of the far
  field for the reference solution equals $2.144$.\label{fig:cube}}
\end{figure}

\begin{figure}[H]
\centering
\includegraphics[height=2.3in]{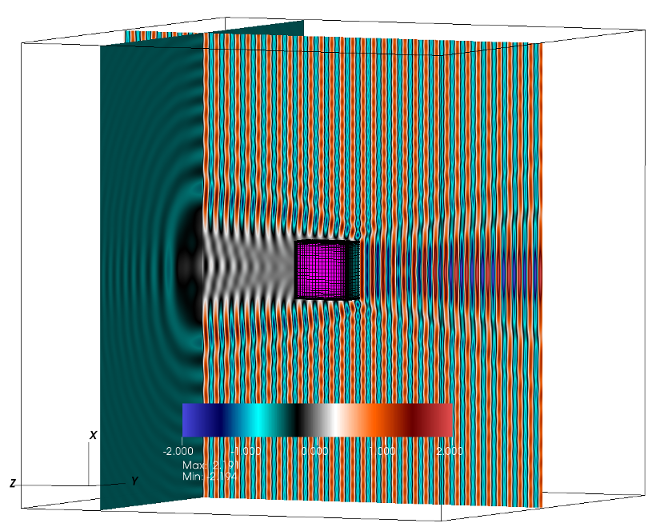} 
\includegraphics[height=2.3in]{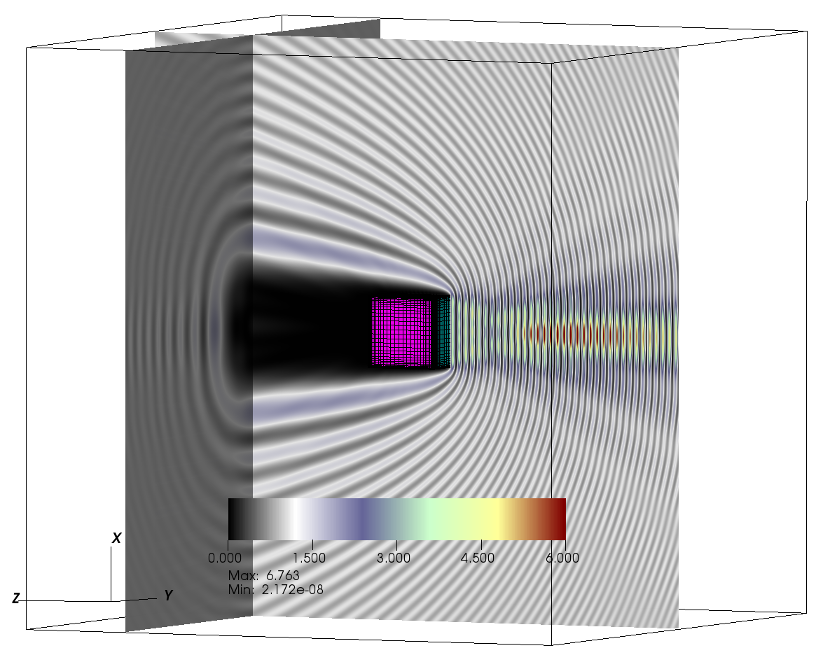} \\
\includegraphics[height=2.3in]{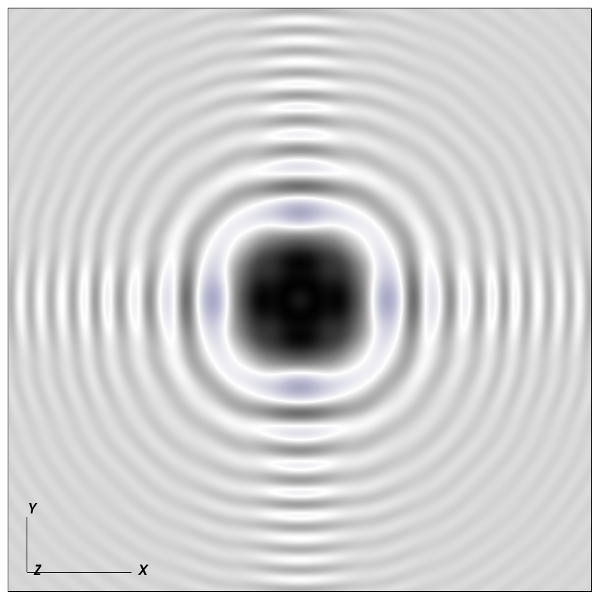} 
\includegraphics[height=2.3in]{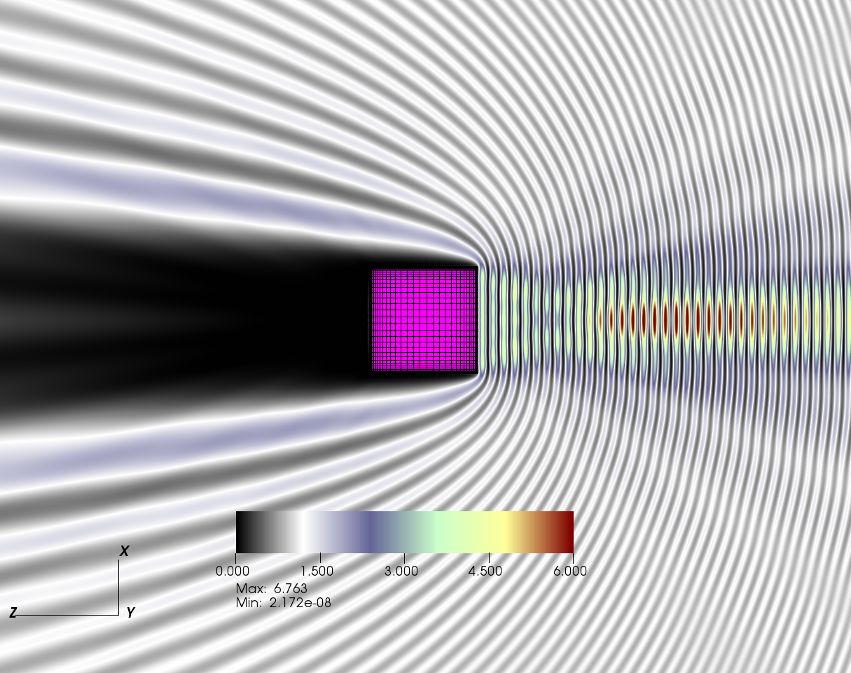} 
\caption{Scattering with $\lambda = 0.4$ by a cube of size $2\times 2 \times 2$
  (each side of size $5\lambda$). The top-left image represents the real part of
  the total field and the other figures show the intensity ($|\field|^2$) from
  different view angles.\label{fig:cube_scat} }
\end{figure}

As mentioned previously, the important problem of scattering by
obstacles containing edges and corners presents several difficulties,
including density and kernel singularities at the edges. In
Figure~\ref{fig:cube}, we demonstrate the performance of the method for
a cube geometry by computing the error in the far field with respect
to a reference solution obtained by using a very fine
discretization. Figure~\ref{fig:cube_scat} shows the scattering
solution by a cube of side $5\lambda$.

\subsection{Open surfaces}
\label{sec:disk}

\begin{figure}[H]
\centering
\includegraphics[scale=0.4]{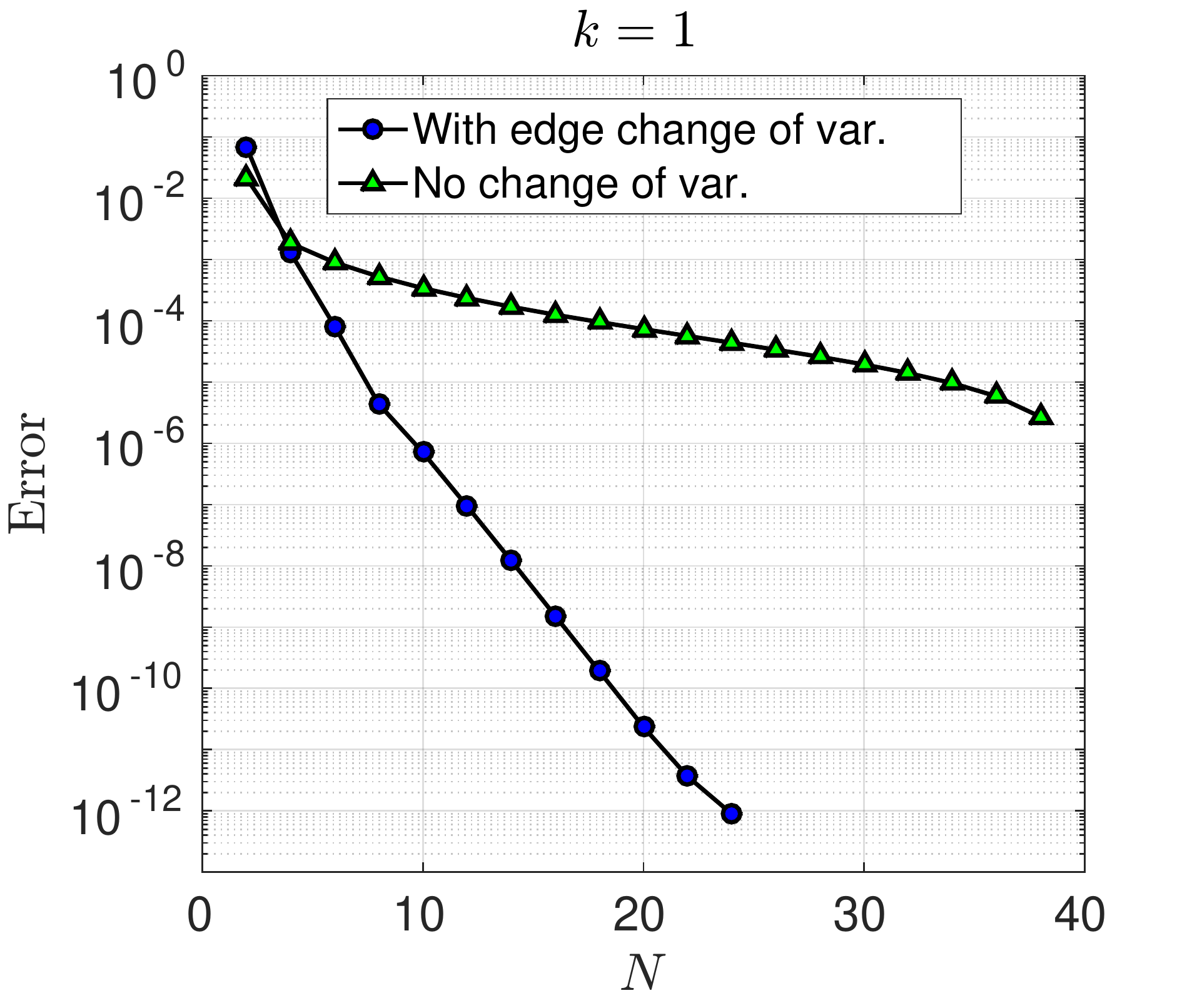}
\caption{Maximum (absolute) far-field error for the problem of scattering by a
  disk of radius 1 and $k=1$. The plot shows both the curve excluding changes of
  variables (in green) and including a $p=4$ change of variables (in blue). The
  maximum value of the far field for the reference solution equals
  $0.7284$.\label{fig:disk}}
\end{figure}

\begin{figure}[H]
\centering
\includegraphics[height=2.3in]{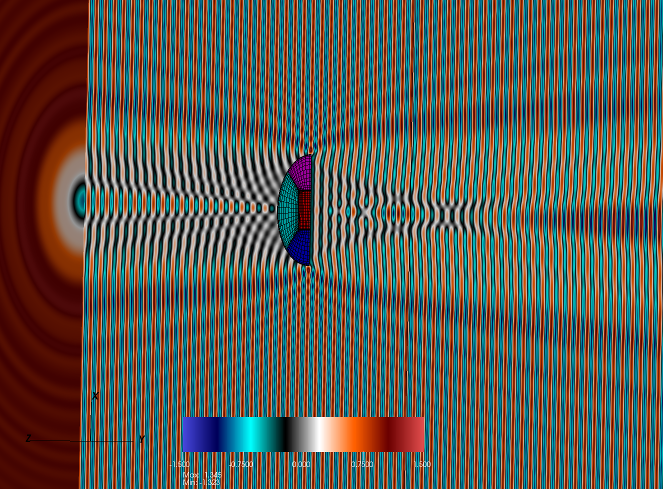} 
\includegraphics[height=2.3in]{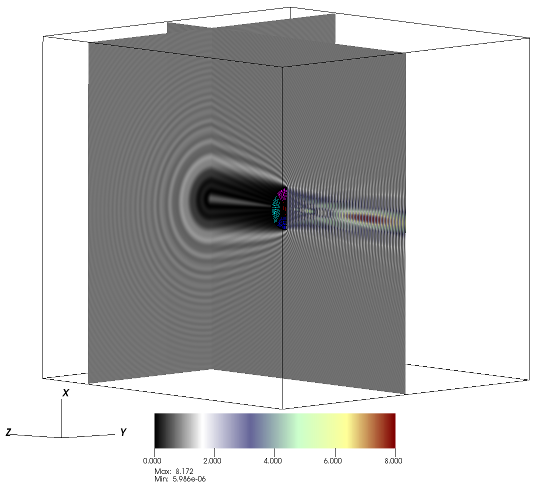} \\
\includegraphics[height=2.3in]{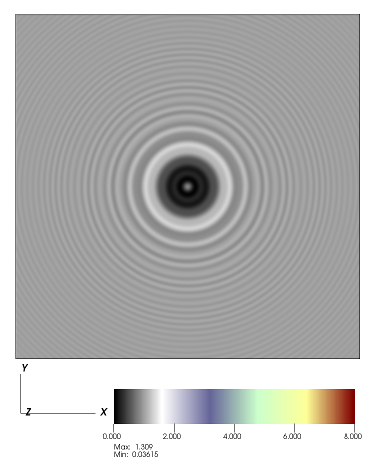} 
\includegraphics[height=2.3in]{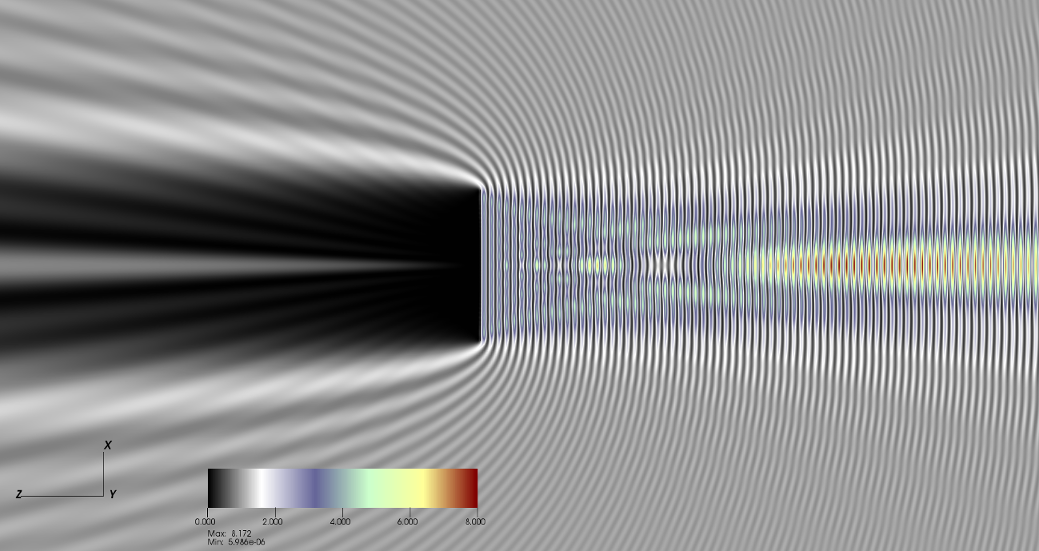} 
\caption{Scattering by a radius-one disk illuminated by a plane wave
  in the $+z$ direction (normal incidence), with $k=10\pi$
  ($\lambda = 0.2$). The image at the top-left displays the real value
  of the total field, while the other images depict the intensity
  ($|\field|^2$) at different view angles. The solution presents the
  well-known Poisson spot (also known as Arago spot and Fresnel bright
  spot) generated at the center of the shaded
  region.\label{fig:disk_scat}}
\end{figure}

Methods for open surfaces typically suffer from low accuracies, or,
alternatively, they require complex treatment at edges. The approach
presented here is a straightforward application of the
rectangular-polar method, with a change of variables at the edges, as
described in Section~\ref{sec:close}. As demonstrated in
Figure~\ref{fig:disk}, which shows the convergence plot for the far
field solution scattered by a disk, the method is robust and
high-order accurate. Figure~\ref{fig:disk_scat} shows the scattering
solution for the problem of scattering by a disk $5\lambda$ in radius.

\subsection{CAD geometries}
\label{sec:cad}

\begin{figure}[H]
\centering
\includegraphics[width=2.0in]{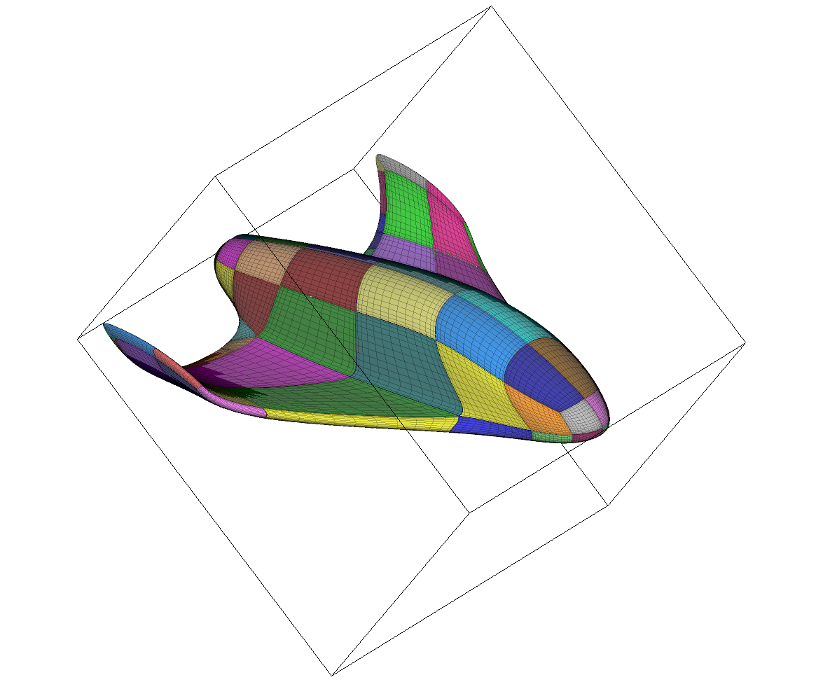} 
\includegraphics[width=2.0in]{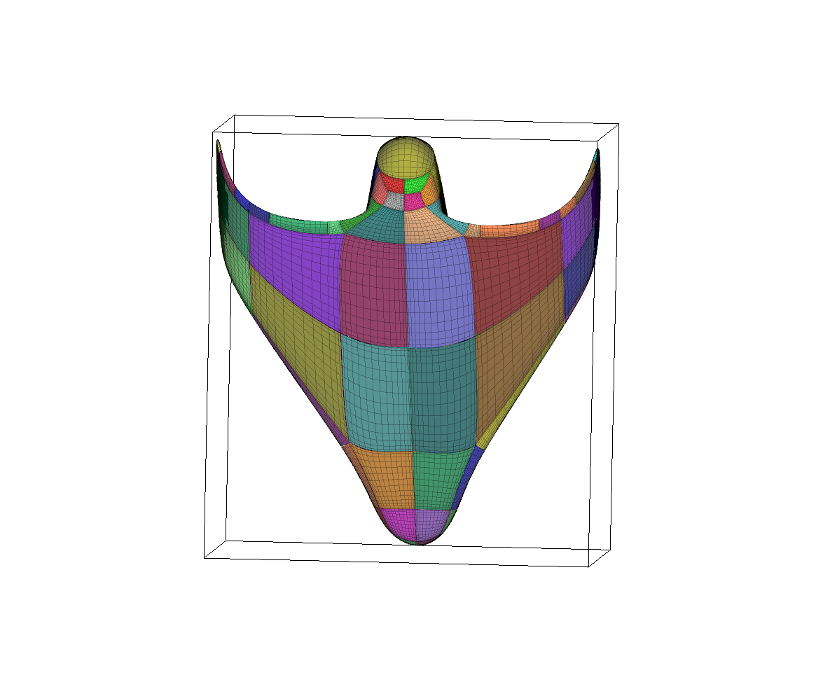} 
\includegraphics[width=2.0in]{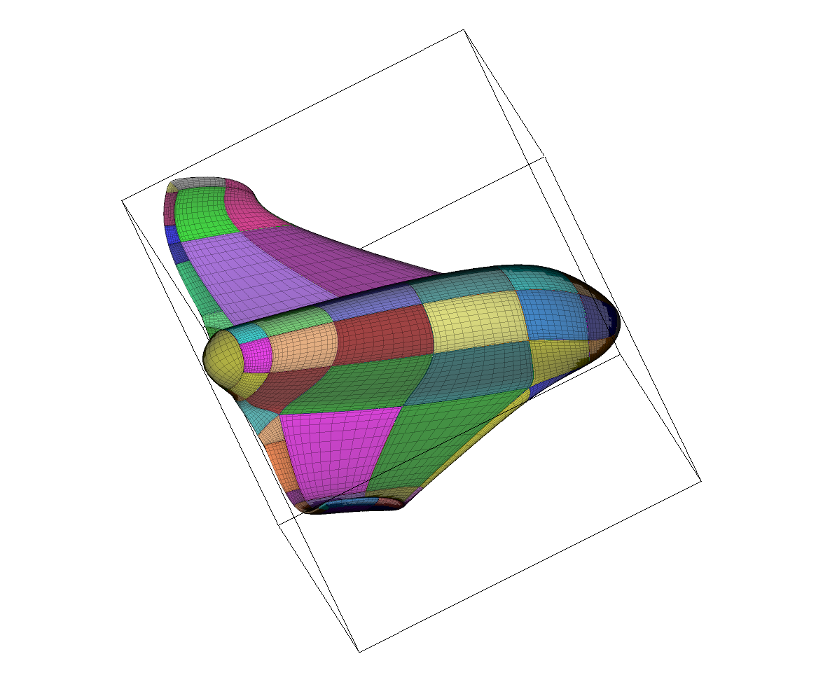}  \\
\includegraphics[width=2.0in]{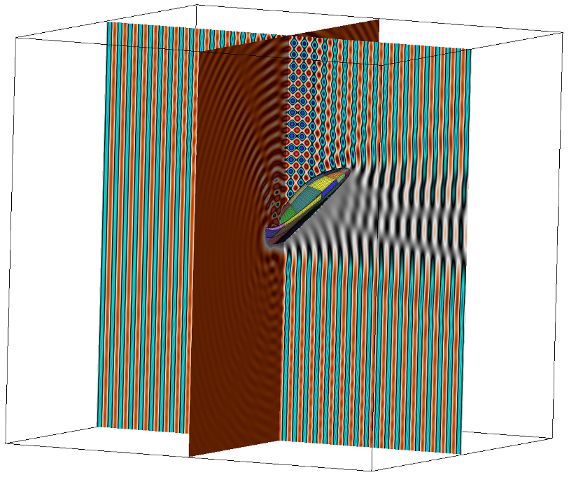} 
\includegraphics[width=2.0in]{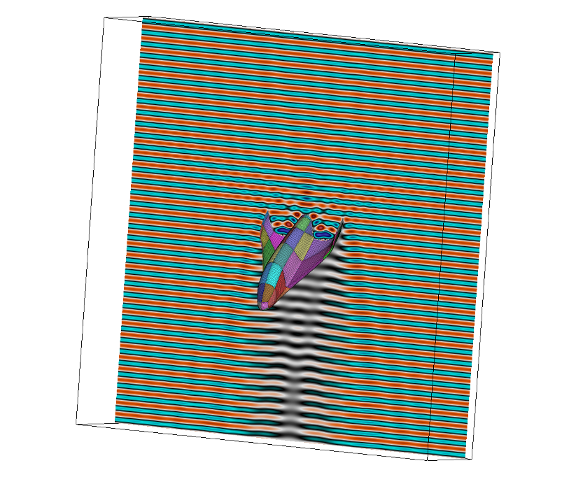} 
\includegraphics[width=2.0in]{glider0007.png} \\
\includegraphics[width=2.0in]{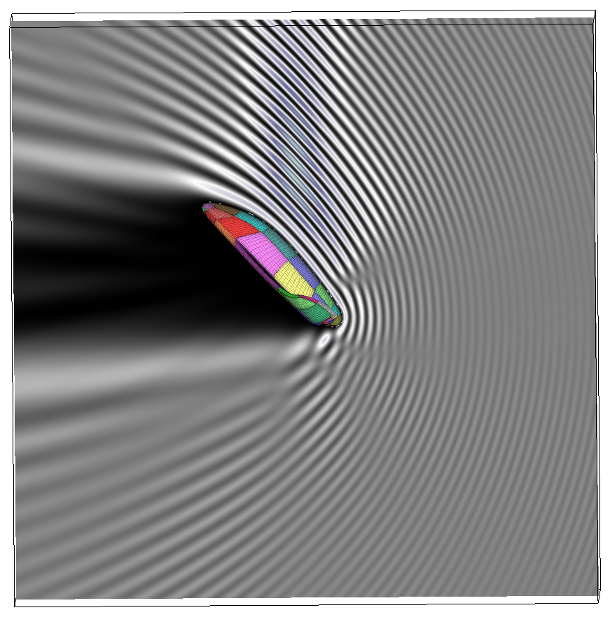} 
\includegraphics[width=2.0in]{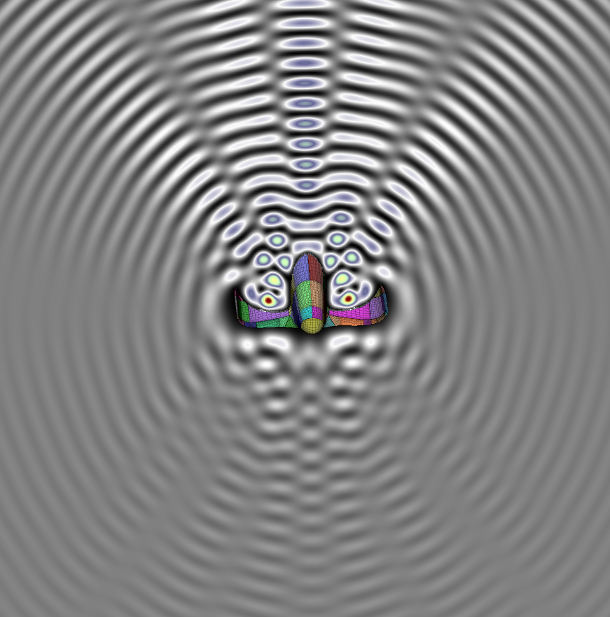} 
\includegraphics[width=2.0in]{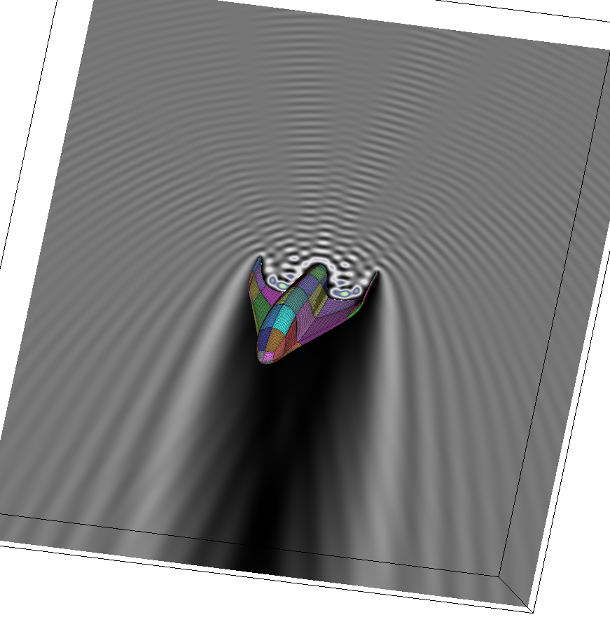}
\caption{Acoustic scattering by a glider geometry from a CAD design. Top row:
  patches and mesh used. Mid row: Real part of the total field. Bottom row:
  Intensity of the total field ($|\field|^2$).\label{fig:glider_scat}}
\end{figure}

As indicated in Section~\ref{surf}, as a rule, CAD designs can easily
be re-expressed as a union of logically-quadrilateral explicitly
parametrized patches, and they are thus particularly well-suited for
use in conjunction with the proposed rectangular-polar solver. To
demonstrate the applicability of the solver to such general type of
geometry descriptions, Figure~\ref{fig:glider_scat} presents the
solution to the problem of acoustic scattering by a glider CAD design
consisting of 87 patches.

\section{Conclusions\label{Conclusions}}

We have presented a rectangular-polar integration strategy for singular
kernels in the context of boundary integral equations. The methodology
was then used in conjunction with the GMRES linear algebra solver to
produce solutions of problems of scattering by obstacles consisting of
open and closed surfaces. As demonstrated by a variety of examples
presented in Section~\ref{numer}, the overall solver produces results
with high-order accuracy, and the rectangular patch description of the
geometry makes the algorithm particularly well-suited for application
to engineering configurations---where the scattering objects are
prescribed in standard (but generally highly complex) CAD
representations.

The proposed methodology has been presented for the context of sound-soft
acoustic scattering, but given the similar nature of the kernels appearing in
the sound-hard case~\citep{Elling-1}, the solver can also be used in that case.
Preliminary results have shown that this solver is also suitable for
electromagnetic scattering in both the perfect electric conductor (PEC) and
dielectric cases.  The numerical examples presented in this paper suggest that
the proposed methodology affords a fast, accurate and versatile high-order
integration and solution methodology for the problem of scattering by arbitrary
engineering structures which, when combined with appropriate acceleration
methods, should result in an accurate solver for highly complex, electrically or
acoustically large problems of propagation and scattering.

\section*{Acknowledgments} %----------------------------------------------------

The authors gratefully acknowledge support by NSF, AFOSR and DARPA
through contracts DMS-1411876 and FA9550-15-1-0043 and HR00111720035,
and the NSSEFF Vannevar Bush Fellowship under contract number
N00014-16-1-2808. Additionally, the authors thank Dr. Edwin Jimenez
and Dr. James Guzman for providing the CAD geometry used in
Section~\ref{sec:cad}.

\bibliographystyle{abbrv}
\bibliography{bibliography}

%-------------------------------------------------------------------------------

\end{document}